%% file: vofini.tex
\documentclass[a4paper,10pt,fleqn]{article}
\usepackage[english]{babel}%

\usepackage[numbers]{natbib}%
\usepackage[bottom=3cm,top=2cm,left=1.5cm,right=1.5cm]{geometry}%
\usepackage{framed,multirow,url,hyperref,setspace,graphicx,bm,xcolor,siunitx,nicefrac,multicol}
\usepackage{amsmath,amssymb,amsthm}
\usepackage{stmaryrd,colortbl}

\usepackage[hang,centerlast]{subfigure}

\usepackage{listings}%

\usepackage{notation}
\usepackage{mathoperators}

\input{vofini_definitions.tex}\usepackage{vofininotation}%
\usepackage{reconstructnotation}%

\let\paragraphold\paragraph%
\renewcommand{\paragraph}[1]{\paragraphold{#1.}}%

\definecolor{newcolor}{rgb}{.8,.349,.1}

\def\worktitle{Third-order accurate initialization of VOF volume fractions on unstructured meshes with arbitrary polyhedral cells}%

\providecommand{\keywords}[1]{\textbf{\textit{Keywords---}} #1}

\allowdisplaybreaks

\begin{document}

\newcommand{\refapp}[1]{appendix~\ref{app:#1}}%

\title{\worktitle}%
\author{Johannes Kromer\href{https://orcid.org/0000-0002-6147-0159}{\includegraphics[height=10pt]{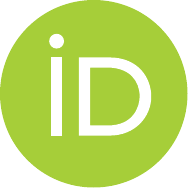}} and Dieter Bothe\textsuperscript{$\dagger$}\href{https://orcid.org/0000-0003-1691-8257}{\includegraphics[height=10pt]{orcid_logo.pdf}}}%
\date{}%
\maketitle
\begin{center}
Mathematical Modeling and Analysis, Technische Universit\"at Darmstadt\\ Alarich-Weiss-Strasse 10, 64287 Darmstadt, Germany\\%
\textsuperscript{$\dagger$}Email for correspondence: \href{mailto:bothe@mma.tu-darmstadt.de?subject=%
Third-order\%20accurate\%20initialization\%20of\%20VOF\%20volume\%20fractions\%20on\%20unstructured%
\%20meshes\%20with\%20arbitrary\%20polyhedral\%20cells}{bothe@mma.tu-darmstadt.de}%
\end{center}

\begin{abstract}%
This paper introduces a novel method for the efficient and accurate computation of volume fractions on unstructured polyhedral meshes, where the phase boundary is an orientable hypersurface, implicitly given as the iso-contour of a sufficiently smooth level-set function. %
Locally, i.e.~in each mesh cell, we compute a principal coordinate system in which the hypersurface can be approximated as the graph of an osculating paraboloid. A recursive application of the \textsc{Gaussian} divergence theorem then allows to analytically transform the volume integrals to curve integrals associated to the polyhedron faces, which can be easily approximated numerically by means of standard \textsc{Gauss-Legendre} quadrature. %
This face-based formulation enables the applicability to unstructured meshes and considerably simplifies the numerical procedure for applications in three spatial dimensions. %
We discuss the theoretical foundations and provide details of the numerical algorithm. Finally, we present numerical results for convex and non-convex hypersurfaces embedded in cuboidal and tetrahedral meshes, showing both high accuracy and third- to fourth-order convergence with spatial resolution. %
\end{abstract}

%
%
\input{01_introduction}%
%
%
\input{02_literature_review}%
%
%
\input{03_mathematical_details}%
%
%
\input{04_numerical_results}%
%
%
\input{05_conclusion}%

%
%
\renewcommand{\bibname}{References}

\input{vofini.bbl}\begin{center}%
\textsc{Acknowledgment}\\[2ex]%
The authors gratefully acknowledge financial support provided by the German Research Foundation (DFG) within the scope of \href{www.sfbtrr75.de}{SFB-TRR 75 (project number 84292822)}. %

The figures in this manuscript were produced using the versatile and powerful library \href{https://ctan.org/topic/pstricks}{\texttt{pstricks}}. For further details and a collection of examples, the reader is referred to the book of \citet{pstricks_2008}.\\[12pt]%
\href{https://www.elsevier.com/authors/policies-and-guidelines/credit-author-statement}{\textsc{CRediT statement}}\\[2ex]%
\textbf{Johannes Kromer}: conceptualization, methodology, software, validation, investigation, data curation, visualisation, writing--original draft preparation, writing--reviewing and editing\\ %
\textbf{Dieter Bothe}: conceptualization, methodology, investigation, writing--reviewing and editing, funding acquisition, project administration%
\end{center}%
%
%
\input{99_appendix}
\end{document}

%% file: vofini_definitions.tex
\newcommand{\reffig}[1]{figure~\ref{fig:#1}}%
\newcommand{\refFig}[1]{Figure~\ref{fig:#1}}%
\newcommand{\reftab}[1]{table~\ref{tab:#1}}%
\newcommand{\refTab}[1]{Table~\ref{tab:#1}}%
\newcommand{\refEqn}[1]{Equation~(\ref{eqn:#1})}%
\newcommand{\refeqs}[1]{eqs.~(\ref{eqn:#1})}%
\newcommand{\refeqno}[1]{(\ref{eqn:#1})}%
\makeatletter
\newcommand{\refeqn}{\@ifstar{\@ifstar\@@refeqn\@refeqn}{\@@@refeqn}}%
\def\@refeqn#1{(\ref{eqn:#1})}%
\def\@@refeqn#1{eqs.~(\ref{eqn:#1})}%
\def\@@@refeqn#1{eq.~(\ref{eqn:#1})}%
\makeatother

\colorlet{shadecolor}{black!12}%
\newtheorem{prototheorem}{Note}[section]
\newenvironment{note}%
{\colorlet{shadecolor}{black!12}\begin{shaded}\begin{prototheorem}}%
{\end{prototheorem}\end{shaded}}%

\newtheorem{remarktheorem}{Remark}[section]%
\newenvironment{remark}%
{\colorlet{shadecolor}{black!12}\begin{shaded}\begin{remarktheorem}}%
{\end{remarktheorem}\end{shaded}}%

\newtheorem{assumptiontheorem}{Assumption}[section]%
\newenvironment{assumption}%
{\colorlet{shadecolor}{black!12}\begin{shaded}\begin{assumptiontheorem}}%
{\end{assumptiontheorem}\end{shaded}}%

\newtheorem*{impltheo}{Implementation detail}
{\colorlet{shadecolor}{magenta!12}\begin{shaded}\begin{impltheo}}%
{\end{impltheo}\end{shaded}}%

%% file: 01_introduction.tex
\section{Introduction}
In the context of a two-phase flow problem in some bounded domain $\domain\subset\setR^3$, the spatial regions $\domain^\pm\fof{t}$ occupied by the respective phases, which are separated by an embedded orientable hypersurface $\iface\fof{t}\subset\domain$, need to be immediately identified. %
One way to achieve this consists in introducing a phase marker $\alpha\fof{t,\vx}$ which, say, is 0 for $\vx\in\domain^+\foft$ and 1 for $\vx\in\domain^-\foft$, respectively. %
A spatial decomposition of the domain into $N_\domain$ pairwise disjoint cells $\cell*_i$ (such that $\domain=\bigcup_i\cell*_i$ with $\cell*_i\cap\cell*_j=\emptyset$ for $i\neq j$) allows to assign to each of those a fraction $\alpha_i:=\lvert\cell*_i\rvert^{-1}\int_{\cell*_i}{\alpha\dvol}$ occupied by the phase $\domain^-$. %
While cells entirely confined in $\domain^-$ and $\domain^+$ exhibit a marker value one and zero, respectively, those intersected by the embedded hypersurface exhibit $0<\alpha_i<1$. %
This representation provides the conceptual foundation of the well-known Volume-of-Fluid (VOF) method introduced by \citet{JCP_1981_vofm}. %
Solving an initial value two-phase flow problem requires, among others, the computation of the aforementioned volume fractions $\alpha_i$ for a given discretized domain $\domain$ and a hypersurface $\iface_{0}$, which describes the initial spatial configuration of the flow. %
If one seeks to compute accurate initial values for curved hypersurfaces, this task becomes particularly challenging, even for seemingly simple hypersurfaces (e.g.\ whose description involves only a small set of parameters) like spheres. Among others, accurate initial values are of paramount importance for the onset of shape instabilities: %
e.g., \citet{JFM_2015_dbob} investigate the dynamic behaviour of high viscosity droplets by releasing initially resting spherical droplets in an ambient liquid. %
Due to buoyancy, the droplets rise and deform, where the rotational symmetry of the configuration quickly degrades for increasing droplet diameter and rise velocity. %
%
Furthermore, accurate volume fractions are required for testing algorithms designed to approximate geometric properties, e.g., curvature and normal fields. %
To the best knowledge of the authors, no higher-order approach applicable to \textit{unstructured} meshes has been published yet. %
In a previous paper \cite{JCP_2019_haco}, the authors have proposed a third-order convergent algorithm for structured  meshes and %
the objective of the present work is to extend the algorithm of \citet{JCP_2019_haco} to unstructured polyhedral meshes. %
Due to the congruence in form and content, the definitions and notation in the remainder of this section as well as the literature review in section~\ref{sec:vof_ini_literature} considerably draw from the respective passages in \cite{JCP_2019_haco}. %
Beyond the extended applicability, the present algorithm also features immanent boundedness (i.e., $0\leq\alpha_i\leq1$) as well as a significant simplification of the numerical procedure. %
We first provide some relevant notation needed to precisely formulate the problem under consideration and to sketch the approach proposed in this work. %
The oriented hypersurface $\iface\subset\domain$ induces a pairwise disjoint decomposition\footnote{Henceforth, we consider a specific instant, say $t=0$, and omit the time argument.} $\domain=\iface\cup\domain^+\cup\domain^-$, where we call $\domain^-$ and $\domain^+$ the \textit{interior} and \textit{exterior} (with respect to $\iface$) subdomain, respectively. %
For the numerical approximation, the domain $\domain$ is decomposed into a set of pairwise disjoint cells $\cell*_i$, some of which are intersected by $\iface$, i.e.\ they contain patches $\iface_i\defeq\iface\cap\cell*_i$ of the hypersurface. %
Any intersected cell again admits a disjoint decomposition into the hypersurface patch $\iface_i$, as well as an "interior" ($\cell*_i^-$) and "exterior" ($\cell*_i^+$) segment. It is important to note that, locally, $\partial\iface_i\neq\emptyset$, even if the hypersurface is globally closed, i.e.\ $\partial\iface=\emptyset$. \refFig{illustration_notation} exemplifies the notation. %
\begin{figure}[htbp]
\null\hfill%
\includegraphics[page=1]{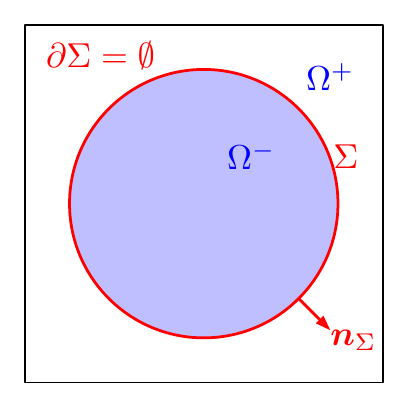}%
\hfill%
\includegraphics[page=2]{domain_decomposition}%
\hfill\null%
\caption{Illustration of the decomposition induced by a closed hypersurface $\iface$.}%
\label{fig:illustration_notation}%
\end{figure}

Henceforth we are concerned with a single intersected polyhedral cell $\cell*_i$ which is why we drop the cell index $i$ for ease of notation. The enclosed patch $\iface$ (i.e., $\iface\cap\cell*_i$) is assumed to be a twice continuously differentiable hypersurface with a piecewise smooth, non-empty boundary $\partial\iface\neq\emptyset$. %
%
%
\begin{note}[Problem formulation]%
For a given polyhedral cell $\cell*$ and hypersurface $\iface$, we seek to approximate %
\begin{align}
\abs{\cell}=\intcell{1}\label{eqn:problem_formulation}%
\end{align}
with high accuracy at finite resolution. %
\end{note}%
%
%
\subsection{Notation}\label{subsec:notation}%
%
%
\paragraph{Computational cells} We consider an arbitrary polyhedron $\cell*$ bounded by $N^\vfface$ planar (possibly non-convex) poly\-gonal faces $\vfface_k$ with outer unit normal $\vn_{\vfface,k}$. To ensure applicability of the \textsc{Gaussian} divergence theorem, $\partial\cell*=\bigcup_k{\vfface_k}$ and $\partial\vfface_k$ are assumed to admit no self-intersections. 
This is not a relevant restriction since objects of the latter class have no relevance for the desired application within finite-volume based methods. %
The $N^\vfface_k$ vertices\footnote{Note that the number of vertices in a closed polygon coincides with the number of edges.} $\set{\vx^\vfface_{k,m}}$ on each face are ordered counter-clockwise with respect to the normal $\vn_{\vfface,k}$, implying that the $m$-th edge $\vfedge_{k,m}$ is spanned by $\vx^\vfface_{k,m}$ and $\vx^\vfface_{k,m+1}$; %
for notational convenience, the indices are continued periodically, i.e.\ $\vx^\vfface_{k,N^{\vfface}_k+1}\defeq\vx^\vfface_{k,1}$. %
\paragraph{Summation} For ease of notation, the summation limits for faces (index $k$) and edges (index $m$) are omitted where no ambiguity can occur. %

%% file: 02_literature_review.tex
%
%
\section{Literature review}\label{sec:vof_ini_literature}%
\begin{figure}[htbp]%
\null\hfill%
\includegraphics[page=1]{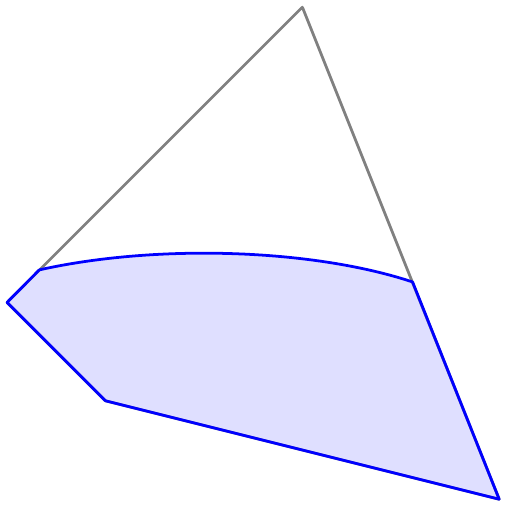}%
\hfill%
\includegraphics[page=2]{lopez_subdivision}%
\hfill%
\includegraphics[page=3]{lopez_subdivision}%
\hfill\null%
\caption{Two-dimensional sketch of the subdivision strategy applied by \citet{JCP_2019_ncaa} with intersected original cell (left), enveloping \textsc{Cartesian} stencil (center, $n_{\mathrm{sc}}=5$) and polyhedral sub-cells composed of hypersurface intersections (\textcolor{blue}{$\bullet$}) and sub-cell vertices (\textcolor{blue}{$\blacksquare$}).}%
\label{fig:lopez_subdivision}%
\end{figure}

The computation of volumes emerging from the intersection of curved hypersurfaces and polyhedral domains (e.g., tetrahedra and hexahedra) has been addressed in several publications up to this date. Besides the direct approaches, i.e.~recursive local grid refinement coupled with a linear approximation of the interface \cite{CS_2005_ecfv,JCP_2010_oric}, the majority of the literature contributions can be classified based on the underlying concept as follows. %
\paragraph{Direct quadrature}%
The work of \citet{CF_2019_tiov} covers the initialization of volume fractions on unstructured grids in two and three spatial dimensions. Their approach consists of decomposing the mesh into simplices and subsequently computing the intersection volume by direct quadrature. The authors report high accuracy for spheres and show that their method is capable of initializing intersections of spheres and hyperboloids, i.e.\ domains with non-smooth boundaries. %
%
%
\citet{IJNME_2016_hoai} develop a higher-order quadrature method for integrals over implicitly defined hypersurfaces, which involves explicitly meshing the zero-isocontour by means of higher-order interface elements. %
\citet{JCP_2016_ecot} propose a computationally efficient and robust method for the computation of volume overlaps of spheres and tetrahedra, wedges and hexahedra. %
The approach of \citet{CF_2015_nioi,CPC_2015_vofi} employs direct computation of integrals with discontinuous integrands by means of quadrature, where the boundaries of the integration domain are computed by a root finding algorithm. While their algorithm requires quite some computational effort, it is capable of handling hypersurfaces with kinks. %
\citet{JCP_2007_gioi} develop an algorithm for geometric integration over irregular domains. To obtain the hypersurface position within an intersected polyhedron, the level-set function is evaluated at the corners, allowing for a linear interpolation along the edges. Subsequent decomposition of the polyhedron into simplices, composed of the interior vertices and intersections, allows for straightforward evaluation of the desired integrals. %
The algorithm of \citet{JCP_2019_ncaa} extends the previous one by a local refinement strategy: each intersected cell is superimposed with a stencil of $n_{\mathrm{sc}}^3$ hexahedral sub-cells, whose respective overlap with the original cell has to be extracted; cf.~\reffig{lopez_subdivision} for an illustration. The authors state that "\textit{taking into account that a higher $n_\mathrm{sc}$ value produces not only a higher initialization accuracy [...] but also a higher CPU time consumed.}" %

\begin{figure}[htbp]
\null\hfill%
\includegraphics[height=2.5cm,page=2]{muller_transformation}\raisebox{1.25cm}{\,$\rightarrow$\,}%
\includegraphics[height=2.5cm,page=4]{muller_transformation}\raisebox{1.25cm}{\,-\,}%
\includegraphics[height=2.5cm,page=3]{muller_transformation}\raisebox{1.25cm}{\,$\rightarrow$\,}%
\includegraphics[height=2.5cm,page=1]{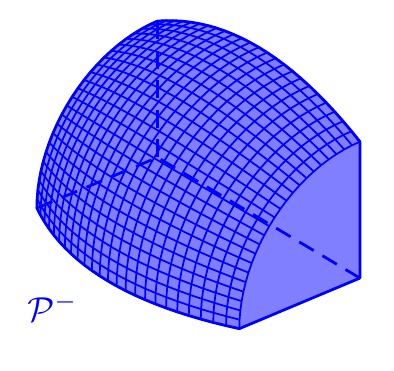}\raisebox{1.25cm}{\,-\,}%
\includegraphics[height=2.5cm,page=3]{muller_transformation}\raisebox{1.25cm}{\,$\rightarrow$\,-}%
\includegraphics[height=2.5cm,page=3]{muller_transformation}%
\hfill\null%
\label{fig:mueller_concept}%
\caption{The domain transformation applied by \citet{IJNME_2013_hasa} exploits that the hypersurface patch $\iface$ can be expressed as the difference of $\partial\protect\cell$ and the union of the immersed faces $\vfface_k^-$. By choosing a set of divergence-free polynomials for the quadrature and applying the \textsc{Gaussian} divergence theorem, the integral over $\protect\cell$ becomes trivial.
}%
\end{figure}

\citet{JCP_2006_tnao} and the series of papers by \citet{JCP_2007_honm,JCP_2009_honm,JSC_2010_honm} are concerned with the numerical evaluation of $\delta$-function integrals in three spatial dimensions. Considering a cuboid intersected by a hypersurface, the concept of Wen is to rewrite the integral over a three-dimensional $\delta$-function as an integral over one of the cell faces, where the integrand is a one-dimensional $\delta$-function. All of the above approaches, however, imply considerable computational effort and complex, case-dependent implementations.
%
\citet{CPC_2005_clfm} introduced a library of four independent routines for multidimensional numerical integration, three of which employ Monte-Carlo integration and the fourth resorts to a globally adaptive subdivision scheme. While methods based on Monte-Carlo integration allow for a wider range of potential applications, the errors exhibit $\hot{N^{-\nicefrac{1}{2}}}$, where $N$ is the number of evaluations of the level-set function, implying comparatively high computational effort to obtain the accuracy desired in most practical applications. %
\paragraph{Divergence theorems}%
\citet{IJNME_2013_hasa} propose an algorithm for the computation of integrals over implicitly given hypersurfaces, resorting to the construction of quadrature nodes and weights from a given level-set function. The computation of a divergence-free basis of polynomials allows reducing the spatial problem dimension by one. By recursive application of this concept, integrals over implicitly defined domains and hypersurfaces in $\setR^3$ are transformed to curve integrals; cf.~\reffig{mueller_concept} for a schematic illustration. While the method of \citet{IJNME_2013_hasa} is computationally highly efficient and exhibits high accuracy, the numerical tests shown by the authors only cover level-set functions of low polynomial order, i.e.\ hypersurfaces with few geometric details and exclusively globally convex ones. Contrary, in section \ref{sec:numerical_results}, we provide results for both locally and globally non-convex hypersurfaces. %
%
%
In a previous contribution \cite{JCP_2019_haco}, we proposed a higher-order method for initialization of volume fractions, based on the combination of a local approximation of the hypersurface by an osculating paraboloid and application of appropriate divergence theorems. The solution of the emerging \textsc{Laplace-Beltrami}-type problem resorts to a \textsc{Petrov-Galerkin} approach, where establishing the linear system of equations requires topological connectivity on a cell level. Beyond this limitation, the algorithm in \cite{JCP_2019_haco} is restricted to simply connected hypersurface patches $\iface\cap\cell*_i$. %
%
%
\paragraph{Discrete hypersurfaces}%
For some applications, such as the breakup of capillary bridges \cite{IJMF_2021_bdoc}, the initial interface configuration results from energy minimization considerations. E.g., the \textit{surface evolver} algorithm of \citet{EM_1992_tse} iteratively approximates the corresponding minimal surfaces by a set of triangles. %
Recently, \citet{XXX_2021_cvfa} proposed an efficient and versatile approach for the initialization of volume fractions on unstructured meshes from such triangulated surfaces. The authors show accurate and second-order convergent results for a variety of triangulated surfaces, including examples with sharp edges and multiple disjoint parts. %
\subsection{Novelties of the proposed approach}%
The novelties of the present approach can be summarized as follows: %
\begin{enumerate}
\item The exploitation of divergence theorems yields an entirely face-based formulation, both supporting efficiency and facilitating the applicability to unstructured meshes with arbitrary polyhedra. %
\item The extended topological admissiblity of the boundary segments $\partial\ifaceapprox$ both eliminates the restriction to simply connected hypersurface patches and allows to handle twice intersected edges, which was not possible in the previous algorithm \cite{JCP_2019_haco}. %
Furthermore, one obtains immanent boundedness, i.e.~$0\leq\alpha_i\leq1$; cf.~\reffig{topological_admissibility}. %
\begin{figure}[htbp]
\null\hfill%
\includegraphics[page=1]{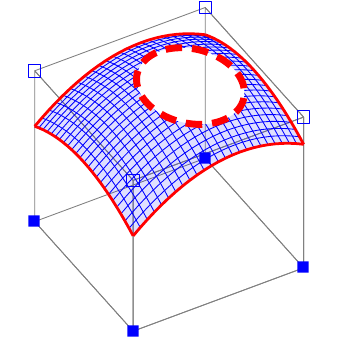}%
\hfill%
\includegraphics[page=2]{clipping_ellipse}%
\hfill%
\includegraphics[page=3]{clipping_ellipse}%
\hfill%
\includegraphics[page=4]{clipping_ellipse}%
\hfill%
\includegraphics[page=5]{clipping_ellipse}%
\hfill\null%
\caption{Novel topological admissibility (dashed) of boundary segments $\partial\ifaceapprox_k$ (left to right): fully enclosed in $\vfface_k$, graph over single edge, disconnected and without graph representation over span of associated intersections (\textcolor{red}{$\bullet$}).}%
\label{fig:topological_admissibility}%
\end{figure}
\item Due to the application of divergence theorems in combination with a local approximation of $\iface$ as the graph of a height function, the proposed method can be easily extended to integrals of type $\int_\iface{\iprod{\vf}{\nS}\darea}$ for functions $\vf$ that are polynomial in the spatial variable $\vx$. %
\end{enumerate}

%% file: 03_mathematical_details.tex
%
%
\section{Mathematical foundations of the method}\label{sec:mathematical_details}%
Let $\cell*\subset\setR^3$ be an arbitrary polyhedron (cf.~subsection~\ref{subsec:notation}), intersected by a twice continuously differentiable oriented hypersurface $\iface$ with outer unit normal $\nS$. The hypersurface is given implicitly as the zero iso-contour of a level-set function $\lvlset_\iface:\setR^3\mapsto\setR$ by %
\begin{align}
\iface=\set{\vx\in\setR^3:\lvlset_\iface\fofx=0}.\label{eqn:hypersurface_levelset}%
\end{align}
For obvious reasons, we assume that $\iface\cap\cell*\neq\emptyset$ and $\iface\not\subset\cell*$, i.e.\ the hypersurface intersects the boundary of the polyhedron $\cell*$. %
We are interested in the volume of the "interior" part of the polyhedron, i.e. %
\begin{align}
\cell\defeq\set{\vx\in\cell*:\lvlset_\iface\fofx\leq0},%
\end{align}
where the superscript "$-$" analogously applies to faces $\vfface_k$ and edges $\vfedge_{k,m}$. %
Henceforth, $\iface$ will be used to abbreviate $\iface\cap\cell*$ for ease of notation. %
Applying the \textsc{Gaussian} divergence theorem and using $\div{\vx}=3$ in $\setR^3$, the volume of an intersected polyhedron $\cell$ (cf.~\refeqn{problem_formulation}) can be cast as
\begin{align}
\abs{\cell}=\intcell{1\,}=\frac{1}{3}\brackets*[s]{\intdcell{\iprod{\vx-\xref}{\vn_{\partial\cell*}}}+\intiface{\iprod{\vx-\xref}{\nS}}},\label{eqn:volume_computation}%
\end{align}
where $\xref$ is an arbitrary but spatially fixed reference point. 
%
With some $\xiref\in\iface$, one may perform a change to the orthonormal base $\set{\vtau_1,\vtau_2,\vn_0}\fof{\xiref}$, where $\vn_0\defeq\nS\fof{\xiref}$ and $\vtau_1$, $\vtau_2$ are the unit normal and principal tangents%
\footnote{I.e., the tangents associated to the principal curvatures $\curv_i$, obtained from the \textsc{Weingarten} map.}%
 to $\iface$ at $\xiref$. For ease of notation, let $\vtau=\brackets[s]{\vtau_1,\vtau_2}$. %
In the vicinity of $\xiref$, the inverse function theorem states that the hypersurface $\iface$ can be expressed as the graph of a height function; see, e.g., the monograph of~\citet{Pruess2016}. In what follows, we assume that $\iface$ admits a unique explicit parametrization as %
\begin{align}
\iface=\set{\vf_\iface\fof{\vt}:\vt\in\support_\iface}\quad\text{with}\quad\vf_\iface\fof{\vt}=\xiref+\vtau\vt+\hiface\fof{\vt}\vn_0,%
\label{eqn:iface_map}%
\end{align}
and some parameter domain $\support_\iface\subset\setR^2$ (henceforth referred to as \textit{graph base} of this parametrization of $\iface$). %
Note that %
\begin{align}
\evaluate{\grad{\hiface}}{\vt=\vec{0}}=\vec{0}%
\quad\text{with}\quad%
\grad{\hiface}=\frac{\partial\hiface}{\partial t_i}\ve_i.%
\end{align}
The tangential coordinates associated to any $\vx\in\setR^3$ are obtained by the projection %
\begin{align}
\vt=\vtau\transpose\brackets{\vx-\xiref},\qquad\text{implying that}\qquad\support_\iface=\set{\vtau\transpose\brackets{\vx-\xiref}:\vx\in\iface}.\label{eqn:tangential_projection}%
\end{align}
%
%
Exploiting $\partial\cell\setminus\iface=\bigcup_{k}\vfface_k^-$ allows to cast the first summand in \refeqn{volume_computation} as %
\begin{align}
\intdcell{\iprod{\vx-\xref}{\vn_{\partial\cell*}}}%
=\sum\limits_{k}{\int_{\vfface_k^-}{\iprod{\vx-\xref}{\vn^\vfface_k}\darea}}%
=\sum\limits_{k}{\iprod{\vx_{k,1}^\vfface-\xref}{\vn^\vfface_k}\immersedarea{k}},\label{eqn:cell_boundary_contribution}%
\end{align}
corresponding to the sum of the immersed face areas $\immersedarea{k}\defeq\abs{\vfface_k^-}$, weighted by the signed distance $\iprod{\vx_{k,1}^\vfface-\xref}{\vn^\vfface_k}$ to the reference $\xref$. %
For the reformulation of the second summand in \refeqn{volume_computation}, note that the explicit parametrization $\iface=\vf_\iface\fof{\support_\iface}$ given in \refeqn{iface_map} allows to express the normal of the hypersurface as %
\begin{align}
\nS%
=\frac{\vn_0-\vtau\grad{\hiface}}{\sqrt{1+\iprod{\grad{}\hiface}{\grad{}\hiface}}}.%
\label{eqn:hypersurface_normal}%
\end{align}
We exploit that $\vx-\xref=\xiref-\xref+\vtau\vt+\hiface\vn_0$ for any $\vx\in\iface$ and \refeqno{hypersurface_normal} to write %
\begin{align}
\iprod{\vx-\xref}{\nS}=%
\frac{\hat{h}_\iface}{\sqrt{1+\iprod{\grad{}\hiface}{\grad{}\hiface}}}%
\quad\text{with}\quad%
\hat{h}_\iface\defeq\hiface-\iprod{\grad{\hiface}}{\vt}+\iprod{\xiref-\xref}{\vn_0-\vtau\grad{\hiface}}.%
\label{eqn:hypersurface_contribution_reduced}%
\end{align}
The integral transformation from $\iface$ to $\support_\iface$ cancels the denomiator in \refeqn{hypersurface_contribution_reduced}, such that one obtains %
\begin{align}
\intiface{\iprod{\vx-\xref}{\nS}}
=\int_{\support_\iface}{\hat{h}_\iface\dd{\vt}}.%
\label{eqn:hypersurface_contribution}%
\end{align}
The continuity of $\hiface$ implies the existence of a function $\Hiface:\setR^2\mapsto\setR^2$ such that $\div{\Hiface}=\hat{h}_\iface$ (a "primitive"). %
Applying the \textsc{Gaussian} divergence theorem once again yields %
\begin{align}
\intiface{\iprod{\vx-\xref}{\nS}}%
=\int_{\support_\iface}{\hat{h}_\iface\dd{\vt}}%
=\int_{\support_\iface}{\div{\Hiface}\dd{\vt}}%
=\int_{\partial\support_\iface}{\iprod{\Hiface}{\vn_{\support_\iface}}\dd{\vt}},%
\label{eqn:hypersurface_contribution_divergence}%
\end{align}
where $\vn_{\support_\iface}$ denotes the outer unit normal to the boundary of the graph base $\support_\iface$. %
Recall that, by assumption, the boundary $\partial\iface=\vf_\iface\fof{\partial\support_\iface}$ of the hypersurface is a subset of the polyhedron boundary, i.e.~$\partial\iface\subset\partial\cell*=\bigcup_k{\vfface_k}$. This suggests a decomposition based on the polyhedron faces $\vfface_k$. Let %
\begin{align}
\partial\iface=\bigcup\limits_{k}{\partial\iface_k}%
\qquad\text{with}\qquad%
\partial\iface_k=\partial\iface\cap\vfface_k=\vf_\iface\fof{\partial\support_{\iface,k}}%
\qquad\text{and}\qquad%
\partial\support_{\iface}=\bigcup\limits_{k}{\partial\support_{\iface,k}},%
\end{align}
which allows to rewrite the rightmost integral in \refeqn{hypersurface_contribution_divergence} as %
\begin{align}
\int_{\partial\support_\iface}{\iprod{\Hiface}{\vn_{\support_\iface}}\dd{\vt}}=%
\sum\limits_{k}{\int_{\partial\support_{\iface,k}}{\iprod{\Hiface}{\vn_{\support_\iface}}\dd{\vt}}}.%
\label{eqn:cell_hypersurface_contribution}%
\end{align}
Finally, combining \refeqs{cell_boundary_contribution} and \refeqno{cell_hypersurface_contribution} yields
\begin{align}
\abs{\cell}=\frac{1}{3}\brackets*[s]{\sum\limits_{k}{\iprod{\vx_{k,1}^\vfface-\xref}{\vn^\vfface_k}\immersedarea{k}+\int_{\partial\support_{\iface,k}}{\iprod{\Hiface}{\vn_{\support_\iface}}\dd{\vt}}}},\label{eqn:volume_computation_facebased}%
\end{align}
implying that the volume of a truncated polyhedral cell $\cell$ can be cast as the sum of face-based quantities (index $k$). %
While the inverse function theorem guarantees the existence of $\hiface$, its actual computation poses a highly non-trival task for general hypersurfaces. However, \refeqs{hypersurface_normal}--\refeqno{cell_hypersurface_contribution} remain valid for an approximated hypersurface 
\begin{align}
\ifaceapprox=\set{\vf_\ifaceapprox\fof{\vt}:\vt\in\support_\ifaceapprox}\quad\text{with}\quad\vf_\ifaceapprox\fof{\vt}=\xiref+\vtau\vt+\happrox\fof{\vt}\vn_0,\label{eqn:ifaceapprox_map}%
\end{align}
where the notation introduced above for $\iface$ analogously applies to $\ifaceapprox$. %
The \textit{principal curvatures} $\curv_i$ and \textit{tangents} $\vtau_i$ at some $\xiref\in\iface$ induce a local second-order approximation %
\begin{align}
\happrox\fof{\vt}=\vt\transpose\vec{\curv}\vt=\frac{\curv_1t_1^2+\curv_2t_2^2}{2}%
=\hiface\fof{\vt}+\hot{\norm{\vt}^3}\qquad\text{with}\qquad%
\curvtensor=\frac{1}{2}\brackets*[s]{\begin{matrix}\curv_1&0\\0&\curv_2\end{matrix}},%
\label{eqn:hypersurface_approximation_quadratic}%
\end{align}
where for $\curv_i\equiv0$ one obtains a tangent plane; subsection~\ref{subsec:approximation_hypersurfaces} describes a procedure to obtain the parameters of the quadratic approximation. %
\begin{remark}[Non-principal approximation]%
The curvature tensor in \refeqn{hypersurface_approximation_quadratic} admits no off-diagonal elements, since it corresponds to a principal coordinate system. However, it is worth noting that the proposed algorithm can be readily extended to the non-principal case, i.e.~with non-zero diagonal elements $\curv_{12}=\curv_{21}$, by replacing \refeqn{hypersurface_approximation_quadratic} with %
\begin{align}
\happrox\fof{\vt}=\vt\transpose\vec{\curv}\vt=\frac{\curv_{11}t_1^2+2\curv_{12}t_1t_2+\curv_{22}t_2^2}{2}%
=\hiface\fof{\vt}+\hot{\norm{\vt}^3}\qquad\text{with}\qquad%
\curvtensor=\frac{1}{2}\brackets*[s]{\begin{matrix}\curv_{11}&\curv_{12}\\\curv_{12}&\curv_{22}\end{matrix}}%
\tag{\ref*{eqn:hypersurface_approximation_quadratic}$^\prime$}.%
\end{align}
\end{remark}%
Note that the hypersurface $\ifaceapprox$ in \refeqn{ifaceapprox_map} can also be expressed implicitly as the zero-isocontour of a level-set, i.e. %
\begin{align}
\ifaceapprox=\set{\vx\in\cell*:\lvlset_\ifaceapprox\fofx=0}%
\quad\text{with}\quad%
\lvlset_\ifaceapprox\fofx=\iprod{\vx-\xiref}{\vn_0}-\brackets{\vx-\xiref}\transpose\vtau\curvtensor\vtau\transpose\brackets{\vx-\xiref}.\label{eqn:ifaceapprox_lvlset}%
\end{align}%
\begin{assumption}[]%
In what follows, we focus on the non-trivial case in which at least one of the principal curvatures (say $\curv_1$) is nonzero. %
\end{assumption}%
For a hypersurface $\ifaceapprox$ of the above class, the integrand in \refeqn{hypersurface_contribution} becomes a third-order polynomial in $\vt$, namely %
\begin{align}
\hat{h}_\ifaceapprox&=h_0-\happrox+\iprod{\vec{\curv}\vec{\xi}}{\vt}\nonumber\\%
&=h_0-\frac{1}{2}\brackets*{\curv_1\brackets{t_1^2+2\xi_1t_1}+\curv_2\brackets{t_2^2+2\xi_2t_2}}%
\quad\text{with}\quad%
h_0\defeq\iprod{\xiref-\xref}{\vn_0}%
\quad\text{and}\quad%
\xi_i=-\iprod{\xiref-\xref}{\vtau_i}.\label{eqn:hypersurface_contribution_coefficients}%
\end{align}
Choosing the reference point to coincide with the paraboloid base point, i.e.~$\xref\defeq\xiref$, implies $h_0=\xi_i=0$. For reasons that will become clear below, we choose the primitive %
\begin{align}
\Happrox\defeq\ve_1\int{\hat{h}_\ifaceapprox\dd{t_1}}%
=-\frac{\ve_1}{6}\hat{H}_\ifaceapprox%
\quad\text{with}\quad%
\hat{H}_\ifaceapprox\defeq\curv_1t_1^3+3\curv_2t_1t_2^2.%
\label{eqn:hypersurface_primitive_paraboloid}%
\end{align}
%
%
%
\begin{note}[The choice of the reference point $\xref$]%
On the one hand, choosing $\xref\defeq\xiref$ apparently implies $h_0=\xi_i=0$ in \refeqn{hypersurface_contribution_coefficients}, implying that the evaluation of \refeqn{hypersurface_primitive_paraboloid} involves fewer multiplications\footnote{The exact gain in efficiency depends, among others, on the compiler options as well the evaluation scheme, e.g., \textsc{Horner}.}. However, $\xiref$ will in general not be coplanar to any of the faces $\vfface_k$, such that the associated immersed areas $\immersedarea{k}$ have to be computed. %
One the other hand, a polyhedral cell $\cell*$ intersected\footnote{Note that, by definition, \textit{intersected} means that the polyhedron $\cell*$ admits at least one interior and one exterior vertex; cf.~\ref{}.} by a paraboloid admits at least three intersected faces. Out of those, at least two, say $\vfface_{k_1}$ and $\vfface_{k_2}$, share a common vertex, implying that the respective containing planes intersect, i.e.~$\abs{\iprod{\vn^\vfface_{k_1}}{\vn^\vfface_{k_2}}}\neq1$. %
While choosing $\xref$ to be coplanar to $\vfface_{k_1}$ and $\vfface_{k_2}$ avoids computing the immersed areas $\immersedarea{k_1}$ and $\immersedarea{k_2}$, the evaluation of \refeqn{hypersurface_primitive_paraboloid} becomes more costly in terms of floating point operations. %
Since the evaluation has to be carried out for all intersected faces, irrespective of the choice of the reference point $\xref$, cf.~\refeqn{volume_computation_facebased}, the most efficient choice can only be substantiated by numerical experiments. %
\end{note}
Replacing the original hypersurface $\iface$ in \refeqn{volume_computation_facebased} by the locally parabolic approximation $\ifaceapprox$ yields the approximative enclosed volume %
\begin{align}
\abs{\cell}\approx\frac{1}{3}\brackets*[s]{\sum\limits_{k}{\iprod{\vx_{k,1}^\vfface-\xref}{\vn^\vfface_k}\immersedarea{k}+%
\int_{\partial\support_{\ifaceapprox,k}}{\iprod{\Happrox}{\vn_{\support_\ifaceapprox}}\dd{\vt}}}}.\label{eqn:volume_computation_facebased_approx}%
\end{align}
Since \refeqn{volume_computation_facebased_approx} expresses the enclosed volume as a sum of face-based quantities, in what follows we focus on a single intersected face $\vfface_k$ of the polyhedron. Before the numerical evaluation \refeqn{volume_computation_facebased_approx} can be addressed, the upcoming subsections discuss the local approximation of hypersurfaces and introduce a classification of the boundary segments $\partial\ifaceapprox_k=\ifaceapprox\cap\vfface_k$. %
%
%
\subsection{A local approximation of the hypersurface $\iface$}\label{subsec:approximation_hypersurfaces}%
In order to obtain a base point $\xiref\in\iface$, on each edge $\vfedge_{k,m}$ we approximate the level-set $\lvlset_\iface$ by a cubic polynomial based on the values and gradients of the level-set evaluated at its respective vertices, i.e.~$\lvlset_\iface\fof{\vx^\vfface_{k,m}}$ and $\grad{\lvlset_\iface}\fof{\vx^\vfface_{k,m}}$. %
\refFig{linear_interpolation_lvlset} illustrates the rationale behind this choice: depending on the sign of the curvature, the linear interpolation of the level-set\footnote{This can be avoided by resorting to a level-set that fulfills the signed distance property. However, finding such a level-set for general hypersurfaces poses a highly non-trivial and thus computationally expensive task in itself.}, as employed by, e.g.,~\citet{JCP_2007_gioi} and \citet{JCP_2019_ncaa}, induces a \textit{systematic} over- or underestimation of the volume fractions that ultimately deteriorates both accuracy and order of convergence. In other words, beyond the planar approximation of a curved hypersurface the linear interpolation of the level-set induces an additional source of volume error. %
\begin{figure}[htbp]%
\null\hfill%
\includegraphics[page=2]{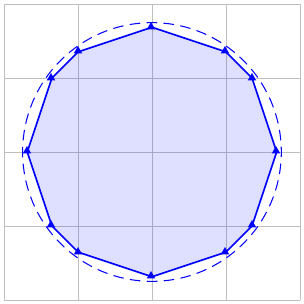}%
\hfill%
\includegraphics[page=3]{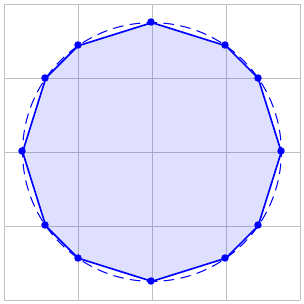}%
\hfill%
\includegraphics[page=4]{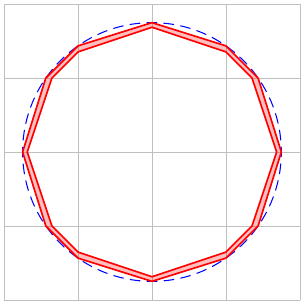}%
\hfill\null%
\caption{Piecewise linear approximation of a circle based on edge intersections obtained from linear (left, \textcolor{blue}{$\blacktriangle$}) and cubic (center, \textcolor{blue}{$\bullet$}) interpolation of the level-set. While the approximation of the circle is piecewise linear in both cases, there is a significant interpolation-induced global area difference (right).}%
\label{fig:linear_interpolation_lvlset}%
\end{figure}%

Note that, in this part of the algorithm, an edge will only be considered \textit{intersected} iff $\lvlset_\iface\fof{\vx^\vfface_{k,m}}\lvlset_\iface\fof{\vx^\vfface_{k,m+1}}<0$, i.e.~if the level-set admits a sign change along the edge $\vfedge_{k,m}$. If existent, the associated root $\vx^\iface_{k,m}$ is computed numerically using a standard \textsc{Newton} scheme. %
\begin{note}[Logical intersection status]%
Technically, deducing the logical intersection status of $\vfface_k$ from the level-set values is not possible for general hypersurfaces. E.g., a face $\vfface_k$ whose vertices $\vx^\vfface_{k,m}$ are entirely located in the negative halfspace of $\iface$ (i.e., $\lvlset_\iface\fof{\vx^\vfface_{k,m+1}}<0$ $\forall m$) may still admit intersections with the positive halfspace of $\iface$. Hence, care has to be taken for the status assignment on the hierarchically superior cell level. In what follows, we assume that the spatial resolution of the underlying mesh is sufficient to capture all geometrical details, implying that each intersected cell contains at least one intersected edge that admits a sign change of the level-set $\lvlset_\iface$. %
\end{note}%
The base point $\xiref\in\iface$ then results from an appropriate projection of the average of all approximate intersections $\set{\vx^\iface_{k,m}}$ onto the hypersurface $\iface$, whose description shall be the content of subsection~\ref{subsec:metric_projection}. %
Finally, the paraboloid parameters in \refeqn{ifaceapprox_map}, namely the normal $\vn_0\defeq\nS\fof{\xiref}$, principal tangents $\vtau\defeq\vtau_\iface\fof{\xiref}$ and curvature tensor $\curvtensor\defeq\curvtensor_\iface\fof{\xiref}$, can be obtained from the \textsc{Weingarten} map. %
\paragraph{Approximation quality}%
The choice of the reference point $\xiref$ crucially affects the global (with respect to the cell) approximation quality of $\ifaceapprox$, measured by the symmetric volume difference (hatched area in \reffig{approximation_quality}). %
\begin{figure}[htbp]%
\null\hfill%
\includegraphics[page=1]{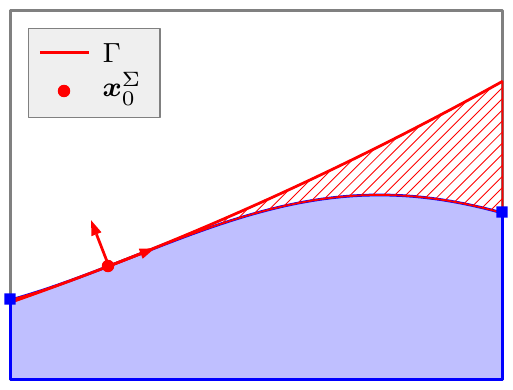}%
\hfill%
\includegraphics[page=3]{approximation_error}%
\hfill%
\includegraphics[page=2]{approximation_error}%
\hfill\null%
\caption{Illustration of qualitative differences in approximation induced by the choice of the base point $\xiref$ (\textcolor{red}{$\bullet$}) in relation to the boundary $\partial\iface$ (\textcolor{blue}{$\blacksquare$}). In the vicinity of the left boundary (left), the quadratic approximation produces a comparatively large symmetric area difference in comparison to a base point with roughly equal distance to the boundary (center). However, there may also be suitable choices of $\xiref$ in the vicinity of the boundary (right).}%
\label{fig:approximation_quality}
\end{figure}%

Under mild restrictions, intuition suggests to select a reference point $\xiref$ close to the (loosely speaking) center of the enclosed hypersurface patch $\iface\cap\cell*$. This aims at reducing the effect of the quadratic growth of the deviation by choosing a reference point whose distance to the boundary $\iface\cap\partial\cell*$ is as uniform as possible. In fact, this provides the motivation behind the projection introduced in subsection~\ref{subsec:metric_projection}, which, as we shall see in section~\ref{sec:numerical_results} below, produces decent results. %
However, as can be seen from the rightmost panel in \reffig{approximation_quality}, this is not necessarily the "best" choice. At this point, note that even the formulation of a minimization problem for general hypersurfaces and polyhedra poses a highly non-trivial task, let alone finding a minimum. Beyond that, due to the cell-wise application, computing such a minimum likely requires considerable computational effort. %
%
%
\subsection{An explicit projection onto the hypersurface $\iface$}\label{subsec:metric_projection}%
The objective of the projection is to compute some point $\xiref\in\iface$ such that the associated paraboloid $\ifaceapprox$ "best" approximates $\iface\cap\cell*$. %
For ease of notation, we (i) relabel the intersections $\vx^\iface_{k,m}$ of the hypersurface $\iface$ with the $m$-th edge $\vfedge_{k,m}$ of the face $\vfface_k$ in a consecutive manner, say $\vy_i$, with $1\leq i\leq I_\iface$ and $I_\iface\defeq\abs{\set{\vx^\iface_{k,m}}}$. %
As illustrated in \reffig{metric_projection}, we employ the following strategy: %
\begin{enumerate}
\item Compute the center as $\bar{\vy}^0=\frac{1}{I_\iface}\sum_{i=1}^{I_\iface}{\vy_i}$ and arrange the shifted $\vy_i$ as the columns of a matrix $\vec{Y}$ of size $3\times I_\iface$ %
by letting $\vec{Y}=\brackets[s]{\vy_1-\bar{\vy}^0,\vy_2-\bar{\vy}^0,\dots,\vy_{I_\iface}-\bar{\vy}^0}$. %
\item Compute the singular value decomposition%
\footnote{In the numerical implementation, we employ the \href{http://www.netlib.org/lapack/explore-html/d1/d7e/group__double_g_esing_ga84fdf22a62b12ff364621e4713ce02f2.html}{LAPACK routine \texttt{DGESVD}}.} (SVD) %
$\vec{Y}=\tensor[2]{U}\tensor[2]{\Sigma}\tensor[2]{V}\transpose$, with the singular values $\tensor[2]{\Sigma}=\diag{\sigma_k}$ such that $\sigma_k\geq\sigma_{k+1}$ and associated left singular vectors $\vu_k=\ve_k\transpose\tensor[2]{U}$. %
After normalization, the left singular vector associated to the \textit{smallest} singular value (i.e.~$\sigma_3$) yields the normal of the plane (containing $\bar{\vy}^0$) that best represents the point cloud $\set{\vy_i}$ in a least-squares sense; hence, let $\vN\defeq\frac{\vu_3}{\sqrt{\iprod{\vu_3}{\vu_3}}}$. %
\item Starting from $\bar{\vy}^0$, we obtain the base point $\xiref$ by iteratively updating %
\begin{align}
\bar{\vy}^{n+1}=\bar{\vy}^n+\vN s^n%
\quad\text{until}\quad%
\abs{\lvlset_\iface\fof{\bar{\vy}^n}}\leq\zerotolerance%
\label{eqn:basepoint_iterative_update}%
,%
\end{align}
where within this work we have chosen $\zerotolerance=\num{e-14}$. %
Projecting the gradient and \textsc{Hessian} of the level-set $\lvlset_\iface$ around $\bar{\vy}^n$ onto the span of $\vN$ yields %
\begin{align}
\lvlset_\iface\fof{\bar{\vy}^n}+\iprod{\grad{\lvlset_\iface}\fof{\bar{\vy}^n}}{\vN}s+\iprod{\vN}{\hessian{\lvlset_\iface}\fof{\bar{\vy}^n}\vN}\frac{s^2}{2}=0.%
\label{eqn:lvl_quad_approx}%
\end{align}
In order to obtain a point on $\iface$ that is close to $\bar{\vy}^0$, the root of \refeqn{lvl_quad_approx} with the smallest absolute value\footnote{Recall from \refeqn{lvl_quad_approx} that $s=0$ corresponds to $\bar{\vy}^n$.} provides the update $s^n$ in \refeqn{basepoint_iterative_update}. With $n^\ast$ denoting the converged iteration, let $\xiref\defeq\bar{\vy}_0^{n^\ast}$. %

\end{enumerate}
\begin{remark}[Principal component analysis (PCA)]%
As an alternative to the singular value decomposition, one could compute the eigenvectors and -values of the $3\times3$ matrix $\vec{Y}\vec{Y}\transpose$. %
To see this, recall that %
$$\vec{Y}\vec{Y}\transpose=%
\tensor[2]{U}\tensor[2]{\Sigma}\tensor[2]{V}\transpose\tensor[2]{V}\tensor[2]{\Sigma}\tensor[2]{U}\transpose=%
\tensor[2]{U}\tensor[2]{\Sigma}^2\tensor[2]{U}\transpose,$$ %
implying that (i) the eigenvalues of $\vec{Y}\vec{Y}\transpose$ are positive and correspond to the square of the singular values $\sigma_k$ and (ii) the eigenvectors correspond to the left singular vectors $\vu_k$. %
Despite being mathematically equivalent, there may be differences in (i) the numerical results due to the accumulation of floating point errors as well as in (ii) the computational time. In the numerical experiments conducted in section~\ref{sec:numerical_results}, we have observed that the singular value decomposition is about 40\% faster than the eigen-decomposition (performed with the \href{http://www.netlib.org/lapack/explore-html/d2/d8a/group__double_s_yeigen_ga442c43fca5493590f8f26cf42fed4044.html\#ga442c43fca5493590f8f26cf42fed4044}{LAPACK routine \texttt{DSYEV}}, including the computation of $\vec{Y}\vec{Y}\transpose$). %
\end{remark}
\begin{figure}[htbp]%
\null\hfill%
\includegraphics{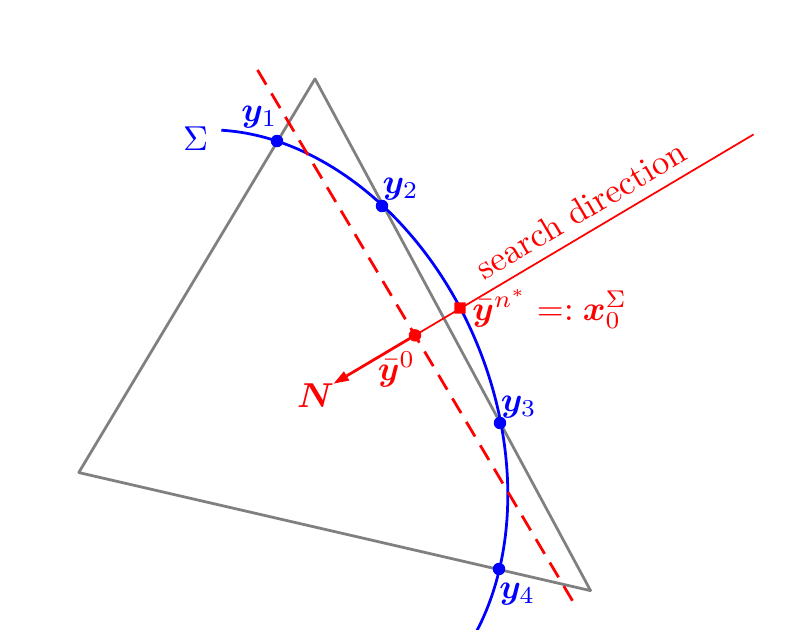}%
\hfill\null%
\caption{Conceptual sketch of the explicit projection of $\bar{\vy}_0$ (\textcolor{red}{$\blacksquare$}) onto the hypersurface $\iface$ (blue).}%
\label{fig:metric_projection}%
\end{figure}%
%
%
\subsection{Intersections of edges $\vfedge_{k,m}$ and paraboloids $\ifaceapprox$}\label{subsec:computation_edge_intersections}%
The intersection of an edge and a paraboloid is conducted in the following way: %
Let the edge $\vfedge_{k,m}$ be parametrized as $$\vfedge_{k,m}=\set*{\vx^\vfface_{k,m}+v\brackets{\vx^\vfface_{k,m+1}-\vx^\vfface_{k,m}}:v\in[0,1]}.$$ %
Due to the quadratic character of the paraboloid $\ifaceapprox$, the projection of the level-set $\lvlset_\ifaceapprox$ onto $\operatorname{span}\fof{\vx^\vfface_{k,m},\vx^\vfface_{k,m+1}}$ can be expressed as a second-order polynomial in the edge coordinate $v$. By inserting the above parametrization into the level-set from \refeqn{ifaceapprox_lvlset}, one obtains %
\begin{align}
\lvlset^\vfedge_{k,m}\fof{v}=&%
\brackets*{\lvlset_\ifaceapprox\fof{\vx^\vfface_{k,m+1}}-\lvlset_\ifaceapprox\fof{\vx^\vfface_{k,m}}-\iprod{\vx^\vfface_{k,m+1}-\vx^\vfface_{k,m}}{\grad{\lvlset_\ifaceapprox}\fof{\vx^\vfface_{k,m}}}}v^2+\nonumber\\%
&\brackets*{\iprod{\vx^\vfface_{k,m+1}-\vx^\vfface_{k,m}}{\grad{\lvlset_\ifaceapprox}\fof{\vx^\vfface_{k,m}}}}v+%
\lvlset_\ifaceapprox\fof{\vx^\vfface_{k,m}}.\label{eqn:ifaceapprox_levelset_edge}%
\end{align}
\refEqn{ifaceapprox_levelset_edge} exhibits $0\leq N^\ifaceapprox_{k,m}\leq2$ real roots $\set{v^\ifaceapprox_{k,m,i}}\subset(0,1)$, to which we associate the intersections $\vx^\ifaceapprox_{k,i}$; cf.~\reffig{intersection_ordering} for an illustration. If two roots are present, we assume without loss of generality that $1>v^\ifaceapprox_{k,m,2}\geq v^\ifaceapprox_{k,m,1}>0$. %
The details of the computation of the relative immersed lengths $\immersedlength{k}{m}$ associated to edge $\vfedge_{k,m}$ can be found in \reftab{immersed_length_edge}. %
\begin{table}[htbp]%
\caption{Relative immersed length $\immersedlength{k}{m}$ (red) on edge $\vfedge_{k,m}$ (gray) intersected with the parabola $\lvlset^\vfedge_{k,m}$; cf.~\protect\refeqn{ifaceapprox_levelset_edge}.}%
\label{tab:immersed_length_edge}%
\centering%
\renewcommand{\arraystretch}{1.5}%
\begin{tabular}{|c||c|c|c|}%
\multicolumn{4}{c}{}\\%
\hline
&\textbf{no intersection}&\textbf{one intersection}&\textbf{two intersections}\\%
&$N^\ifaceapprox_{k,m}=0$&$N^\ifaceapprox_{k,m}=1$&$N^\ifaceapprox_{k,m}=2$\\%
\hline\hline%
$\lvlset^\vfedge_{k,m}\fof{0}<0$&$0$&$1-v^\ifaceapprox_{k,m,1}$&$v^\ifaceapprox_{k,m,2}-v^\ifaceapprox_{k,m,1}$\\%
&\includegraphics[page=4]{roots_on_edge}&\includegraphics[page=5]{roots_on_edge}&\includegraphics[page=1]{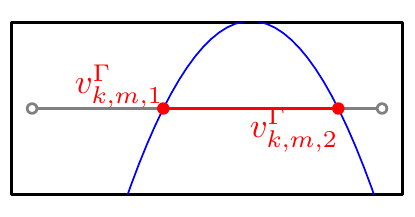}\\%
\hline%
$\lvlset^\vfedge_{k,m}\fof{0}>0$&$1$&$v^\ifaceapprox_{k,m,1}$&$1-v^\ifaceapprox_{k,m,2}+v^\ifaceapprox_{k,m,1}$\\%
&\includegraphics[page=3]{roots_on_edge}&\includegraphics[page=6]{roots_on_edge}&\includegraphics[page=2]{roots_on_edge}\\%
\hline%
\end{tabular}%
\end{table}%

Each intersection $\vx^\ifaceapprox_{k,i}$ can be classified based on the sign of $\evaluate{\partial_v\lvlset^\vfedge_{k,m}}{v=v_{0,i}}$ as either \textit{entering} (negative) or \textit{leaving} (positive), providing a key information for establishing the topological connectivity required for the computation of immersed areas $\immersedarea{k}$ in subsection~\ref{subsec:topological_connectivity}. %
Note that the total number of intersections on each face $\vfface_k$ is even\footnote{Intersected vertices must hence be either considered twice or not all all.} and corresponds to twice the number of curved segments $M_k$, i.e.~$2M_k=\sum_{m=1}^{N^\vfface_k}{N^\ifaceapprox_{k,m}}$. %
Furthermore, the sequence of intersections $\set{\vx^\ifaceapprox_{k,m}}$ alternates between \textit{entering} and \textit{leaving} if the edges $\vfedge_{k,m}$ of the face are traversed in counter-clockwise order with respect to the face normal $\vn^\vfface_k$. %
%
%
\subsection{Transformation to principal coordinates}\label{subsec:transformation_principal_coordinates}%
The robust treatment of the curved segments $\partial\ifaceapprox_k$ calls for a classification of the intersection curves in terms of locally principal coordinates; cf.~\reftab{possible_curve_class}. %
Let the polyhedron face $\vfface_k$ be parametrized as %
\begin{align}
\vfface_k=\set{\vx^\vfface_{k,1}+\vmu_k\vu:\vu\in\support_{\vfface,k}}%
\quad\text{with}\quad%
\vmu_k=\brackets[s]{\vmu_{k,1},\vmu_{k,2}},\quad%
\vmu_k\transpose\vmu_k=\ident,\quad%
\vmu_k\transpose\vn^\vfface_k=\vec{0},%
\label{eqn:cell_face_parametrization}%
\end{align}
and some parameter domain $\support_{\vfface,k}\subset\setR^2$. %
We apply the tangential projection from \refeqn{tangential_projection} to the map given in \refeqn{cell_face_parametrization} and plug the result into \refeqn{ifaceapprox_map} to obtain the implicit quadratic definition of the boundary curve segment $\partial\ifaceapprox_k=\ifaceapprox\cap\vfface_k$, namely %
\begin{align}
\partial\ifaceapprox_k=\set{\vu\in\support_{\vfface,k}:\iprod{\vu}{\tA_k\vu}+\iprod{\vu}{\va_k}+a_k=0}%
\label{eqn:support_levelset_definition}%
\end{align}
with the coefficients %
\begin{align}
a_k=\iprod{\vx^\vfface_{k,1}-\xiref}{\vtau\curvtensor\vtau\transpose\brackets{\vx^\vfface_{k,1}-\xiref}-\vn_0},\quad%
\va_k=\vmu_k\transpose\brackets*{2\vtau\curvtensor\vtau\transpose\brackets{\vx^\vfface_{k,1}-\xiref}-\vn_0}\quad\text{and}\quad%
\tA_k=\vmu_k\transpose\vtau\curvtensor\vtau\transpose\vmu_k.\!\!%
\end{align}
\refTab{possible_curve_class} gathers and illustrates the admissible curve classes that emerge from \refeqn{support_levelset_definition}. %
\begin{table}[htbp]%
\centering%
\renewcommand{\arraystretch}{1.25}%
\caption{Classification of intersections $\partial\ifaceapprox_k$ (red) of planar faces $\vfface_k$ (red shaded) and paraboloids $\ifaceapprox$ (blue). Note that the linear case requires the tangential plane and the face to intersect orthogonally.}%
\label{tab:possible_curve_class}%
\begin{tabular}{l|c|c|c|c}
\multicolumn{5}{c}{}\\%
&\textbf{hyperbolic}&\textbf{elliptic}&\textbf{parabolic}&\textbf{linear}\\%
&\includegraphics[width=3cm,page=1]{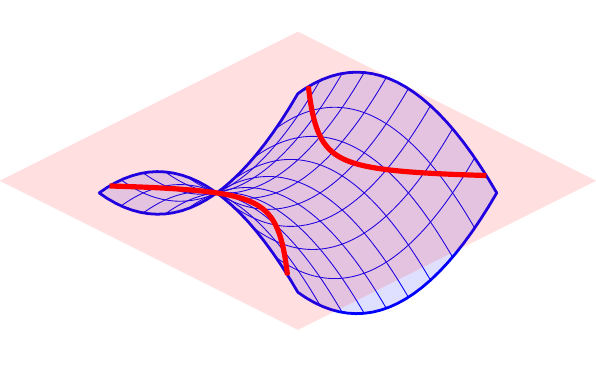}%
&\includegraphics[width=3cm,page=3]{curve_classes}%
&\includegraphics[width=3cm,page=2]{curve_classes}%
&\includegraphics[width=3cm,page=4]{curve_classes}\\%
\hline%
$\textcolor{black}{0<}\abs{\iprod{\vn^\vfface_k}{\vn_0}}\leq1$&\cellcolor{black!25}$\curv_1\curv_2<0$&%
\cellcolor{black!25}$\curv_1\curv_2>0$&\cellcolor{black!25}$\begin{matrix}\curv_1\curv_2=0\\\curv_1+\curv_2\neq0\end{matrix}$&\\%
\hline%
$\phantom{0<}\abs{\iprod{\vn^\vfface_k}{\vn_0}}=0$&&&\cellcolor{black!25}$\curv_1\curv_2\lesseqqgtr0$&\cellcolor{black!25}$\curv_1\curv_2=0$%
\end{tabular}
\end{table}

Note that the matrix of quadratic coefficients $\tA_k$ will not admit diagonal form in general. %
With%
\footnote{The proof of this statement employs that, by definition, we have $\norm{\vmu_{k,1}}=\norm{\vmu_{k,2}}=1$, $\iprod{\vmu_{k,1}}{\vmu_{k,2}}=0$ and $\iprod{\vtau_1}{\vtau_2}=0$. %
Expanding the $i$-th row of $\vtau\transpose\vmu_k=\vec{0}$ yields $\iprod{\vtau_i}{\vmu_1}=\iprod{\vtau_i}{\vmu_2}=0$, implying that %
$\vtau_i=\beta_i\brackets{\cross{\vmu_1}{\vmu_2}}=\beta_i\vn^\vfface_k$ with some $\beta_i\in\setR\setminus\set{0}$. In combination, one obtains the contradiction $\iprod{\vtau_1}{\vtau_2}=\beta_1\beta_2\neq0$.} %
$\vtau\transpose\vmu_k\neq\vec{0}$, the principal coordinates $\vu^\prime$ emerge from $\vu$ via %
\begin{align}
\vu=\tensor[2]{R}_k\transpose\vu^\prime+\vu_{0,k}%
\quad\text{with}\quad%
\tensor[2]{R}_k=\brackets*[s]{\begin{matrix}\cos\psi_k&\sin\psi_k\\-\sin\psi_k&\cos\psi_k\end{matrix}}%
\quad\text{and}\quad%
\tan2\psi_k=\frac{2A_{k,12}}{A_{k,11}-A_{k,22}}.%
\label{eqn:principal_face_coordinates}%
\end{align}
For notational convenience, let $\vr_{i,k}=\tensor[2]{R}_k\transpose\ve_i$. %
With \refeqn{principal_face_coordinates}, the quadratic equation in \refeqn{support_levelset_definition} can be rewritten as %
\begin{align}
\iprod{\vu^\prime}{\tensor[2]{B}_k\vu^\prime}+\iprod{\vu^\prime}{\vb_k}+b_k=0%
\label{eqn:support_levelset_definition_principal}%
\end{align}
with the coefficients  %
\begin{align}
b_k=a_k+\iprod{\vu_{0,k}}{\va_k+\tA_k\vu_{0,k}},%
\quad%
\vb_k=\tensor[2]{R}_k\brackets{\va_k+2\tA_k\vu_{0,k}}%
\quad\text{and}\quad%
\tensor[2]{B}_k=\tensor[2]{R}_k\tA_k\tensor[2]{R}_k\transpose=\diag{B_{k,1},B_{k,2}}.%
\label{eqn:support_levelset_definition_principal_coeff}%
\end{align}
The eigenvalues of $\tensor[2]{B}_k$ ($B_{k,1}\geq B_{k,2}$) classify the boundary segment $\partial\ifaceapprox_k$: %
\begin{enumerate}
\item For $\det\tensor[2]{B}_k\neq0$, the curve segment is \textit{elliptic} ($\det\tensor[2]{B}_k>0$) or \textit{hyperbolic} ($\det\tensor[2]{B}_k<0$). In both cases, the coefficients in \refeqn{support_levelset_definition_principal_coeff} read %
\begin{align}
\vb_k=\vec{0}%
\quad\text{and}\quad%
b_k=a_k+\frac{\iprod{\va_k}{\vu_{0,k}}}{2}%
\quad\text{with the origin}\quad%
\vu_{0,k}=-\frac{1}{2}\tA_k\inverse\va_k.%
\label{eqn:curveclass_det_nonzero}%
\end{align}
%
%
\item For $\det{\tensor[2]{B}_k}=\det{\tA_k}=0$, the curve segment is either \textit{parabolic} ($B_{k,1}\neq0$, $B_{k,2}=0$) or \textit{linear} ($B_{k,1}=B_{k,2}=0$). %
For parabolic intersections, exploiting that $\iprod{\vr_{2,k}}{\tA_k\vr_{2,k}}=B_{k,2}=0$ yields 
\begin{align}
\vb_k=\tensor[2]{R}_k\va_k%
\quad\text{and}\quad%
b_k=0%
\quad\text{with the origin}\quad%
\vu_{0,k}=%
\begin{cases}%
-\frac{a_k}{\iprod{\va_k}{\vr_{2,k}}}\vr_{2,k}&\text{if}\quad\iprod{\va_k}{\vr_{2,k}}\neq0,\\%
-\frac{\iprod{\va_k}{\vr_{1,k}}}{2B_{k,1}}\vr_{1,k}&\text{if}\quad\iprod{\va_k}{\vr_{2,k}}=0.%
\end{cases}%
\label{eqn:curveclass_det_zero}%
\end{align}
In the second case of \refeqn{curveclass_det_zero}, the intersection consists of two parallel lines. While this corresponds to a parabola whose vertex is located at infinity, we prefer a treatment as a degenerate hyperbola for consistency of implementation. %
For linear intersections ($B_{k,1}=B_{k,2}=0$), one obtains %
\begin{align}%
\vb_k=\tensor[2]{R}_k\va_k%
\quad\text{and}\quad%
b_k=a_k%
\quad\text{with the origin}\quad%
\vu_{0,k}=\vec{0}.%
\end{align}%
\end{enumerate}
%
%
With the curve parameter $s$, the above classification induces the following explicit parametrizations: %
\begin{align}
\vu_k^\prime\fof{s}=%
\begin{cases}
\brackets*[s]{\sqrt{\frac{-B_{k,1}}{b_k}}\cos s,\sqrt{\frac{-B_{k,2}}{b_k}}\sin s}\transpose&\text{elliptic},\\%
\brackets*[s]{\pm\sqrt{-\frac{b_k+B_{k,2}s^2}{B_{k,1}}},s}\transpose&\text{hyperbolic},\\%
\brackets*[s]{s,\frac{\iprod{\va_k}{\vr_{1,k}}}{\iprod{\va_k}{\vr_{2,k}}}s+\frac{B_{k,1}}{\iprod{\va_k}{\vr_{2,k}}}s^2}\transpose&\text{parabolic},\\%
\brackets*[s]{s\sin\varphi_k-b_k\cos\varphi_k,-s\cos\varphi_k-b_k\sin\varphi_k}\transpose&\text{linear with }\tan\varphi_k=\frac{\iprod{\va_k}{\vr_{2,k}}}{\iprod{\va_k}{\vr_{1,k}}}.%
\end{cases}
\label{eqn:parametrization_curve_class}%
\end{align}
By plugging \refeqn{parametrization_curve_class} into \refeqn{principal_face_coordinates}, one obtains an explicit parametrization of the boundary segment, i.e. %
\begin{align}%
\partial\ifaceapprox_k=%
\set*{\vx^\vfface_{k,1}+\tensor[2]{\mu}_k\brackets*{\vu_{k,0}+\tensor[2]{R}_k\vu_k^\prime\fof{s}}:s\in\support_{\partial\ifaceapprox,k}}%
\quad\text{with}\quad%
\support_{\partial\ifaceapprox,k}=\bigcup_{m=1}^{M_k}\brackets[s]{s_{k,2m-1}^\ifaceapprox,s_{k,2m}^\ifaceapprox},%
\label{eqn:boundary_curve_principal_coordinates}%
\end{align}
where the union over $M_k$ (number of curved segments of $\partial\vfface_k^-$) intervals reflects the fact that $\partial\ifaceapprox_k$ is not necessarily simply connected; cf.~the rightmost panel in \reffig{face_intersection_elliptic_hyperbolic}. %
The interval boundaries $s_{k,m}^\ifaceapprox$ are obtained by first projecting the edge intersections $\vx^\ifaceapprox_{k,m}$ (cf.~subsection \ref{subsec:computation_edge_intersections}) onto the principal coordinates of the face $\vfface_k$ via %
\begin{align}
\vu^\prime_{k,m}=\tensor[2]{R}_k\brackets*{\tensor[2]{\mu}_k\transpose\brackets*{\vx^\ifaceapprox_{k,m}-\vx^\vfface_{k,1}}-\vu_{0,k}}%
\end{align}
and subsequent inversion of the respective parametrization in \refeqn{parametrization_curve_class}. %
%
%
Projecting the map in \refeqn{boundary_curve_principal_coordinates} onto the base plane of the paraboloid $\ifaceapprox$ using \refeqn{tangential_projection} yields %
\begin{align}
\partial\support_{\ifaceapprox,k}=%
\set{\vt_{0,k}+\tensor[2]{T}_{0,k}\vu_k^\prime\fof{s}:s\in\support_{\partial\ifaceapprox,k}}%
\quad\text{with}\quad%
\vt_{0,k}\defeq\vtau\transpose\brackets{\vx^\vfface_{k,1}-\xiref+\tensor[2]{\mu}_k\vu_{k,0}}%
\quad\text{and}\quad%
\tensor[2]{T}_{0,k}\defeq\vtau\transpose\tensor[2]{\mu}_k\tensor[2]{R}_k,%
\label{eqn:boundary_support_segment}%
\end{align}
corresponding to the integration domain required for the evaluation of \refeqn{volume_computation_facebased_approx}. %
From \refeqn{boundary_support_segment}, the boundary normals are obtained via %
\begin{align}
\vn_{\support_{\ifaceapprox,k}}=%
\brackets*[s]{\ve_2,-\ve_1}%
\frac{\tensor[2]{T}_{0,k}\partial_s\vu_k^\prime}%
{\sqrt{\iprod{\tensor[2]{T}_{0,k}\partial_s\vu_k^\prime}{\tensor[2]{T}_{0,k}\partial_s\vu_k^\prime}}},%
\quad\text{where}\quad%
\iprod*{\vn^\vfface_k-\frac{\iprod{\vn^\vfface_k}{\grad{\lvlset_\ifaceapprox}}}{\iprod{\grad{\lvlset_\ifaceapprox}}{\grad{\lvlset_\ifaceapprox}}}\grad{\lvlset_\ifaceapprox}}{\vtau\vn_{\support_\ifaceapprox,k}}\stackrel{!}{>}0%
\label{eqn:boundary_support_outer_normal}%
\end{align}
must be enforced by sign inversion (if needed) to ensure that $\vn_{\support_{\ifaceapprox,k}}$ is an outer normal to $\support_\ifaceapprox$; cf.~\reffig{illustration_graphbase} for an illustration. %
Recall from \refeqn{hypersurface_primitive_paraboloid} that, by design, $\iprod{\Happrox}{\ve_2}\equiv0$. This allows to rewrite the second summand in \refeqn{volume_computation_facebased_approx} as %
\begin{align}
\int\limits_{\partial\support_{\ifaceapprox,k}}{\iprod{\Happrox}{\vn_{\support_\ifaceapprox}}\dd{\vt}}=%
-\sum\limits_{m=1}^{M_k}{\int\limits_{s_{k,2m-1}^\ifaceapprox}^{s_{k,2m}^\ifaceapprox}{\ve_2\transpose\tensor[2]{T}_{0,k}\partial_s\vu_k^\prime\fof{s}\hat{H}_\ifaceapprox\fof{\vt_k\fof{s}}\dd{s}}}%
\quad\text{with}\quad%
\vt_k\fof{s}=\vt_{0,k}+\tensor[2]{T}_{0,k}\vu_k^\prime\fof{s},%
\label{eqn:quadrature_line_integral}%
\end{align}
which will be evaluated by standard \textsc{Gauss-Legendre} quadrature using $N_{\mathrm{quad}}=5$ nodes. %
\begin{figure}[htbp]
\null\hfill%
\includegraphics{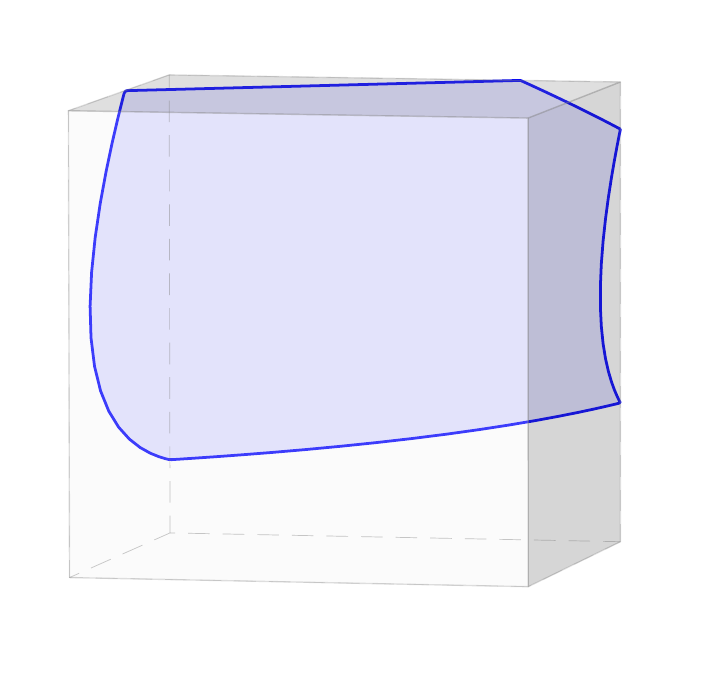}%
\hfill%
\includegraphics{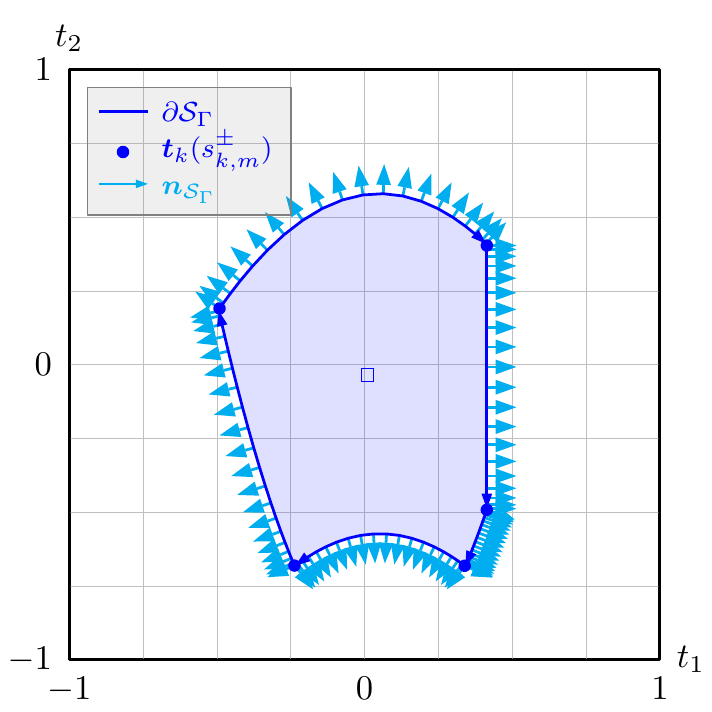}%
\hfill\null%
\caption{Illustration of intersected hexahedral cell with associated graph base $\support_\ifaceapprox$.}%
\label{fig:illustration_graphbase}%
\end{figure}

\subsection{Topological connectivity of curve segments}\label{subsec:topological_connectivity}%
\label{subsec:topological_connectivity}%
The logical status of the principal origin $\vu_{k,0}$, denoted $\vfstatus{}_k$, induces the correct orientation for elliptic, hyperbolic and linear curve segments. In the parabolic case, the origin $\vu_{k,0}$ is located on the curve segment. Hence, one needs to consider the focus point of the parabola to assess the orientation ($\textcolor{blue}{\blacksquare}$ in the bottom left panel in \reffig{intersection_ordering}). One obtains %
\begin{align}%
\vfstatus{}_k=%
\begin{cases}%
-\sign{b_k}&\text{elliptic/hyperbolic/parabolic (degenerate)},\\%
-\sign{a_k}&\text{linear},\\%
\phantom{-}\sign{B_{k,1}}&\text{parabolic}.%
\end{cases}
\label{eqn:logical_status_center}%
\end{align}%
For elliptic\footnote{A proper implementation of the arcus tangens ensures that the transition of the numerical values for the angles at $0$ and $2\pi$ is handled correctly.}, parabolic and linear segments on \textbf{convex} faces $\vfface_k$, traversing the edges $\vfedge_{k,m}$ counter-clockwise with respect to the normal $\vn^\vfface_k$ yields a properly ordered sequence of intersections $\set{\vx^\ifaceapprox_{k,m}}$. %
Hyperbolic segments additionaly require to first assign the intersections to the respective branch of the hyperbola (cyan/blue in the top left panel of \reffig{intersection_ordering}). For an exterior center $\vu_{0,k}$ ($\vfstatus{}_k=+1$ in \refeqn{logical_status_center} and $\square$ in \reffig{intersection_ordering}), the order of the sequence center must be inverted. %
After the intersections have been arranged in this manner, one needs to ensure that the first intersection is of type \textit{entering} by possibly performing an index shift\footnote{Note that the direction of the shift is irrelevant for our purpose. However, we shift to the left, i.e.~$\vx^\ifaceapprox_{k,m}\mapsto\vx^\ifaceapprox_{k,m-1}$ for all $m$.}. \refFig{intersection_ordering} illustrates the concept. %
\begin{figure}[htb]%
\null\hfill%
\includegraphics[page=5]{boundary_curve_ordering}%
\hfill%
\includegraphics[page=3]{boundary_curve_ordering}%
\hfill\null%
\\%
\null\hfill%
\includegraphics[page=1]{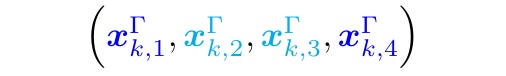}%
\hfill%
\includegraphics[page=5]{boundary_curve_ordering_annotation}%
\hfill\null%
\\%
\null\hfill%
\includegraphics[page=2]{boundary_curve_ordering_annotation}%
\hfill%
\includegraphics[page=6]{boundary_curve_ordering_annotation}%
\hfill\null%
\\%
\null\hfill%
\includegraphics[page=1]{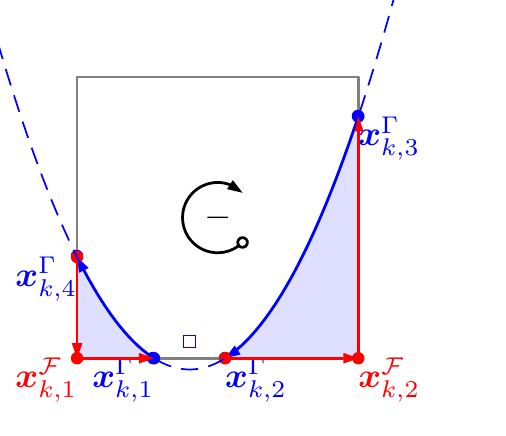}%
\hfill%
\includegraphics{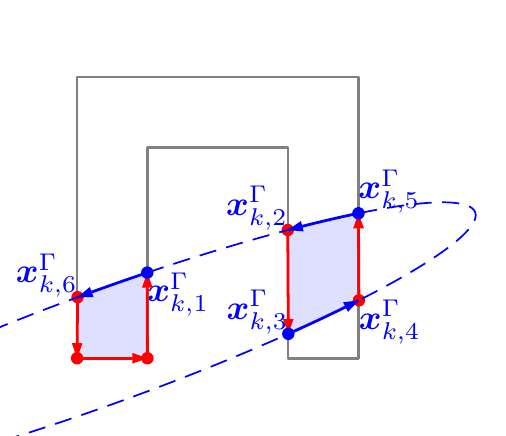}%
\hfill\null%
\\%
\null\hfill%
\includegraphics[page=3]{boundary_curve_ordering_annotation}%
\hfill%
\includegraphics[page=1]{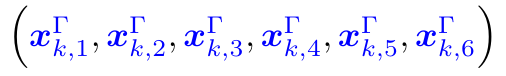}%
\hfill\null%
\\%
\null\hfill%
\includegraphics[page=4]{boundary_curve_ordering_annotation}%
\hfill%
\includegraphics[page=2]{boundary_curve_ordering_nonconvex_annotation}%
\hfill\null%
\caption{Intersection ordering for different curve classes (hyperbolic/convex, elliptic/convex, parabolic/convex and elliptic/non-convex) with initial (top) and sorted (bottom) list of intersections. Immersed edges are colored red whereas boundary curve segments are colored blue (blue/cyan for the branches of the hyperbola).}%
\label{fig:intersection_ordering}%
\end{figure}
\begin{remark}[Non-convex faces $\vfface_k$]%
As can be seen from the bottom right panel in \reffig{intersection_ordering}, non-convex faces potentially degrade the immanent order of intersections. %
However, the correct order of the intersections can be established be arranging the curve parameters $s^\pm_{k,m}$ associated to $\vx^\ifaceapprox_{k,m}$ (cf.~\refeqn{parametrization_curve_class}) in an ascending/descending order. %
\end{remark}%
For the remainder of this manuscript, we assume that the intersections $\set{\vx^\ifaceapprox_{k,m}}$ are arranged in the above manner. %
\subsection{Computation of immersed areas $\immersedarea{k}$}\label{subsec:computation_immersed_area}%
In a manner similar to \refeqn{volume_computation}, the boundary of an immersed face $\vfface_k^-$ can be decomposed into linear and curved segments, where each of the latter connects two edge intersections. %
For each curved segment, we introduce an edge $\vfedge^\ifaceapprox_{k,m}$ (green in the center panel of \reffig{face_decomposition}) that connects the associated intersections. In other words, we (i) replace the curved segments by lines to form a (set of) polygons and (ii) connect the end points of the removed curved segments to form closed paths confining curved "caps". %
\begin{figure}[htbp]%
\null\hfill%
\includegraphics[page=1]{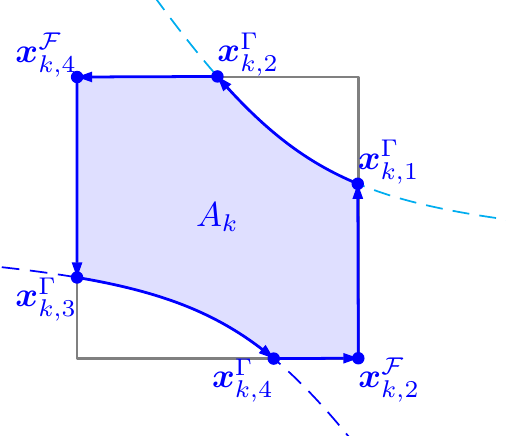}%
\hfill%
\includegraphics[page=2]{boundary_face_decomposition}%
\hfill%
\includegraphics[page=3]{boundary_face_decomposition}%
\hfill\null%
\caption{Decomposition of immersed face $\vfface_k^-$ (left) into polygonal (center) and non-polygonal (right) part.}%
\label{fig:face_decomposition}%
\end{figure}%

This implies that the immersed area $\immersedarea{k}$ is fed from two contributions: %
\begin{enumerate}%
\item The \textbf{polygonal} part can be computed from two sets of edges: (i) the immersed segments of the original edges $\vfedge^-_{k,m}$, which we represent by the original edge $\vfedge_{k,m}$ and its associated relative immersed length $\immersedlength{k}{m}$ (cf.~\reftab{immersed_length_edge}) and (ii) the $M_k$ curve segment bases $\set{\vfedge^\ifaceapprox_{k,m}}=\set{\brackets{\vx^\ifaceapprox_{k,2m-1},\vx^\ifaceapprox_{k,2m}}}$ introduced before, whose arrangement was described in subsection~\ref{subsec:topological_connectivity}. %
\item The \textbf{non-polygonal} part can be computed efficiently by resorting to the principal transformation introduced in subsection~\ref{subsec:transformation_principal_coordinates}%
; cf.~\reffig{face_decomposition} for an illustration. %
\end{enumerate}%
Note that, while the polygonal contribution $\immersedarea{k}^{\mathrm{poly}}$ is zero or positive, the non-polygonal contribution $\immersedarea{k}^{\mathrm{cap}}$ may become negative if the immersed face $\vfface_k^-$ is non-convex (as shown in the right panel of \reffig{face_decomposition}). %
%
%
\paragraph{Polygonal segment ($\immersedarea{k}^{\mathrm{poly}}$)}%
Applying the \textsc{Gaussian} divergence theorem to a polygon embedded in $\vfface_k$ yields %
\begin{align}
\immersedarea{k}^{\mathrm{poly}}=\frac{1}{2}\brackets*{%
\sum\limits_{m=1}^{N^\vfface_k}{\iprod{\vx_{k,m}^\vfface}{\cross{\brackets{\vx_{k,m+1}^\vfface-\vx_{k,m}^\vfface}}{\vn^\vfface_k}}\immersedlength{k}{m}}+%
\sum\limits_{m=1}^{M_k}{\iprod{\vx_{k,2m-1}^\ifaceapprox}{\cross{\brackets{\vx_{k,2m}^\ifaceapprox-\vx_{k,2m-1}^\ifaceapprox}}{\vn^\vfface_k}}}%
},%
\label{eqn:polygonal_area_computation}%
\end{align}
with the relative immersed lengths $\immersedlength{k}{m}$ from \reftab{immersed_length_edge}. %
However, recall that the computation of the area of a planar polygon in $\setR^3$ actually poses a two-dimensional problem. %
Following \citet{JCP_2016_anvc,JCP_2019_ncaa,CPC_2020_voft} and \citet{JCP_2021_fbip}, who for their part resort to the work of \citet{JGT_2002_fpaa}, we employ a projection onto one of the coordinate planes (with normal $\vn_p$ and $\iprod{\vn^\vfface_k}{\vn_p}\neq0$): %
\begin{align}
\abs{\vfface_k}=\frac{\abs{\vfface^{\mathrm{proj}}_{k,p}}}{\iprod{\vn^\vfface_k}{\vn_p}}%
\qquad\text{with}\qquad%
p\defeq\arg\max_p\abs{\iprod{\vn^\vfface_k}{\ve_p}};\label{eqn:face_area_projection}%
\end{align}
see~\reffig{coordinate_projection} for an illustration. %
In order for the projection to maintain the counter-clockwise order of the vertices (i.e., with respect to $\vn_p$), the projected coordinates must be arrangend as $\set{y,z}$ for $p=1$, $\set{z,x}$ for $p=2$ and $\set{x,y}$ for $p=3$. %
\begin{figure}[htbp]%
\null\hfill%
\includegraphics[page=1]{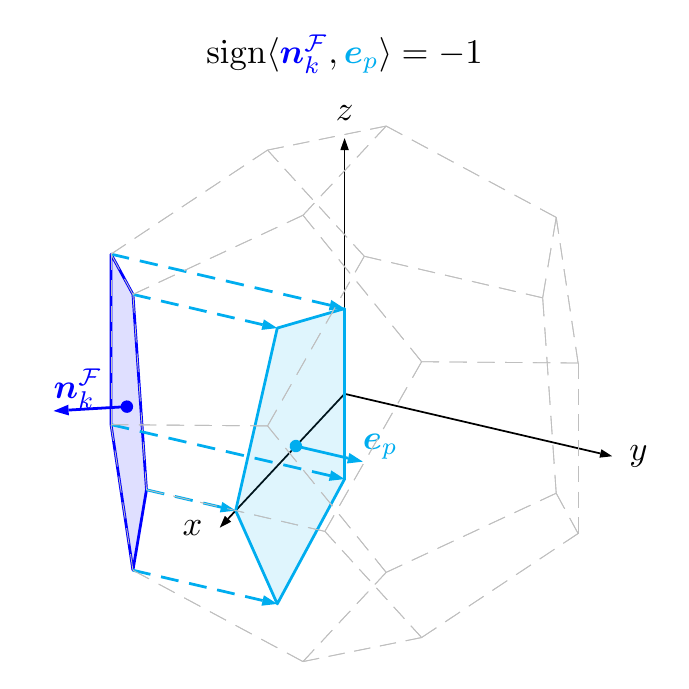}%
\hfill%
\includegraphics[page=2]{projected_face_polygon}%
\hfill\null%
\caption{Coordinate plane projection $\vfface^{\mathrm{proj}}_{k,p}$ (cyan) of polygonal face $\vfface_k$ (blue). Note that the coordinate arrangement produces positive areas in \protect\refeqn{polygonal_area_computation} by accounting for $\operatorname{sign}\protect\iprod{\vn^\vfface_k}{\vn_p}$.}%
\label{fig:coordinate_projection}%
\end{figure}
Since the projection acts on $\darea$ as a scalar multiplication, one can simply substitute the three-dimensional points in \refeqn{polygonal_area_computation} with their projected counterparts, where the exterior vector product reduces to swapping two vector components. %
E.g., for a projection onto the $xy$-plane ($p=3$), \refeqn{polygonal_area_computation} becomes %
\begin{align}
\immersedarea{k}^{\mathrm{poly}}=\frac{1}{2\iprod{\vn^\vfface_k}{\ve_3}}%
\brackets*{%
\sum\limits_{m=1}^{N^\vfface_k}{\brackets*{x^\vfface_{k,m}y^\vfface_{k,m+1}-y^\vfface_{k,m}x^\vfface_{k,m+1}}%
\immersedlength{k}{m}}%
+%
\sum\limits_{m=1}^{M_k}{\brackets*{x^\ifaceapprox_{k,2m-1}y^\ifaceapprox_{k,2m}-y^\ifaceapprox_{k,2m-1}x^\ifaceapprox_{k,2m}}}%
},%
\label{eqn:polygonal_area_computation_projected}%
\end{align}
and analogous expressions for the $xz$- and $yz$-plane. %
%
%
\paragraph{Non-polygonal segment ($\immersedarea{k}^{\mathrm{cap}}$)}%
The parametrization of the curve segment in principal coordinates $\vu^\prime_k\fof{s}$ can be obtained from \refeqn{parametrization_curve_class}. %
We obtain %
\begin{align}%
\immersedarea{k}^{\mathrm{cap}}=%
\begin{cases}%
\frac{\sqrt{B_{k,1}B_{k,2}}}{2b_k}\sum\limits_{m=1}^{M_k}{s_{k,2m}^\ifaceapprox-s_{k,2m-1}^\ifaceapprox + \cos{s_{k,2m}^\ifaceapprox}\sin{s_{k,2m-1}^\ifaceapprox}-\cos{s_{k,2m-1}^\ifaceapprox}\sin{s_{k,2m}^\ifaceapprox}}&\text{elliptic},\\%
\sqrt{\frac{b_k}{4B_{k,1}}}\sum\limits_{m=1}^{M_k}{s_{k,2m-1}^\ifaceapprox\lambda_{2m}-s_{k,2m}^\ifaceapprox\lambda_{2m-1}+%
\sqrt{\frac{-b_k}{B_{k,2}}}\brackets{\rho_{2m}-\rho_{2m-1}}}&\text{hyperbolic ($+$)},\\%
\sqrt{\frac{b_k}{4B_{k,1}}}\sum\limits_{m=1}^{M_k}{s_{k,2m}^\ifaceapprox\lambda_{2m-1}-s_{k,2m-1}^\ifaceapprox\lambda_{2m}-%
\sqrt{\frac{-b_k}{B_{k,2}}}\brackets{\rho_{2m}-\rho_{2m-1}}}&\text{hyperbolic ($-$)},\\%
\frac{B_{k,1}}{6\iprod{\va_k}{\vr_{2,k}}}\sum\limits_{m=1}^{M_k}{\brackets{s_{k,2m}^\ifaceapprox-s_{k,2m-1}^\ifaceapprox}^3}&\text{parabolic},\\%
0&\text{linear}%
\end{cases}%
\label{eqn:cap_area}%
\end{align}
with $\lambda_m\defeq\sqrt{1-\frac{B_{k,2}}{b_k}\brackets*{s_{k,m}^\ifaceapprox}^2}$ and %
$\rho_{m}\defeq\log\brackets*{\sqrt{-\frac{B_{k,2}}{b_k}}\lambda_m-\frac{B_{k,2}}{b_k}s_{k,m}^\ifaceapprox}$ for ease of notation; %
cf.~\refFig{face_intersection_elliptic_hyperbolic} for an illustration. %
\begin{figure}[htbp]%
\null\hfill%
\includegraphics[page=2]{curve_class_illustration_face}%
\hfill%
\includegraphics[page=1]{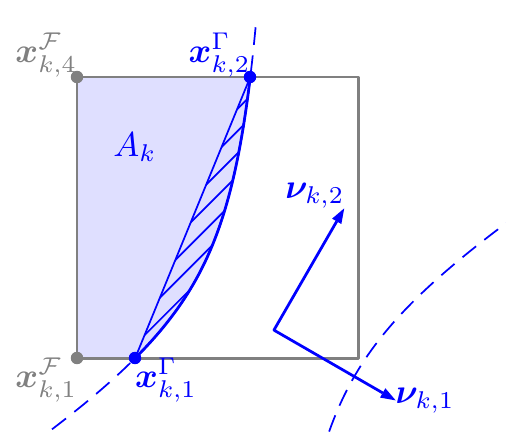}%
\hfill%
\includegraphics[page=3]{curve_class_illustration_face}%
\hfill\null
\caption{Elliptic, hyperbolic and parabolic face intersection ($\vec{\nu}_k=\tensor[2]{R}_k\tensor[2]{\mu}_k$ denote the principal coordinates associated to the face $\vfface_k$) with immersed area $\immersedarea{k}$ and cap area $\immersedarea{k}^{\mathrm{cap}}$ (hatched) from \protect\refeqn{cap_area}.}%
\label{fig:face_intersection_elliptic_hyperbolic}%
\end{figure}
\paragraph{Fully enclosed boundary segments}%
For faces $\vfface_k$ containing hyperbolic and parabolic boundary curve segments $\partial\ifaceapprox_k$, the absence of edge intersections implies either $\immersedarea{k}=0$ if all vertices are exterior or $\immersedarea{k}=\abs{\vfface_k}$ if all vertices are interior. %
Contrarily, as can be seen from \reffig{elliptic_interiority}, an ellipse, say $\mathcal{E}$, can either (i) be fully enclosed in the face ($\mathcal{E}\subset\vfface_k$, left panel), (ii) fully enclose the face ($\mathcal{E}\supset\vfface_k$, center panel) or (iii) admit no overlap ($\mathcal{E}\cap\vfface_k=\emptyset$, right panel). Note that the classification, which is of paramount importance for the topological admissibly (cf.~\reffig{topological_admissibility}), cannot be deduced from the status of the vertices but only from the status of the center $\vx^\vfface_{k,1}+\tensor[2]{\mu}_k\vu_{k,0}$; cf.~\refeqn{parametrization_curve_class}. %

\begin{figure}[htbp]%
\null\hfill%
\includegraphics{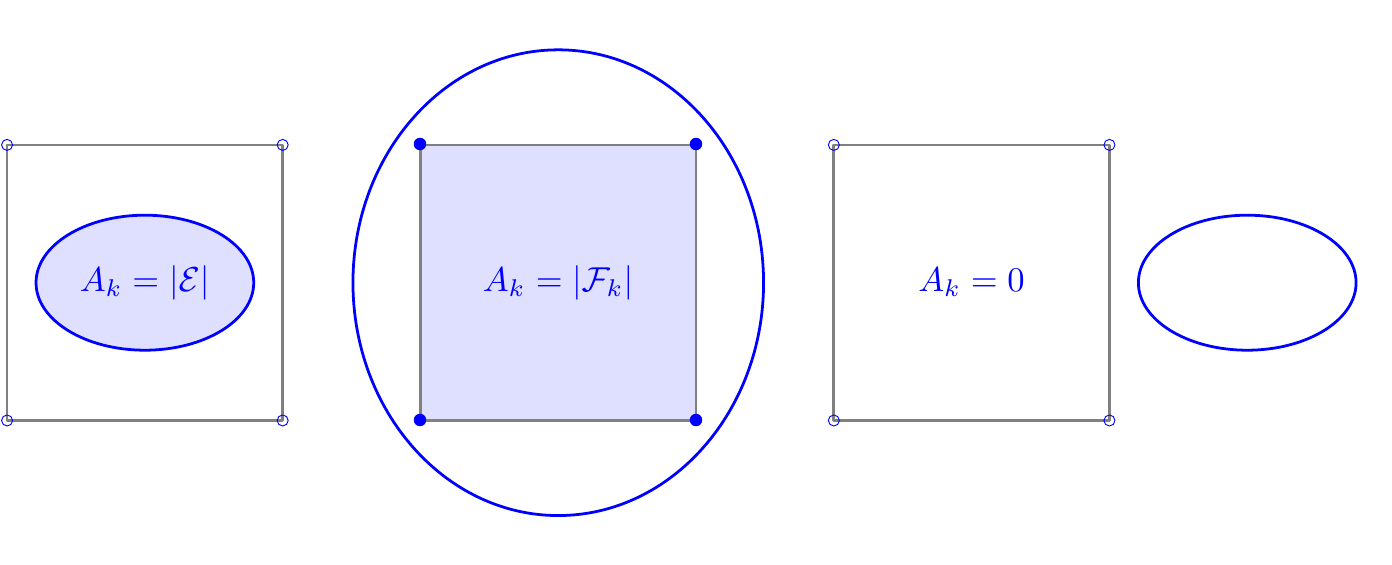}%
\hfill\null%
\caption{Ellipse $\mathcal{E}$ that is fully enclosed, fully encloses and admits no overlap, where interior (exterior) vertices of the containing face $\vfface_k$ are denoted by $\textcolor{blue}{\bullet}$ ($\textcolor{blue}{\circ}$); cf.~\protect\reffig{topological_admissibility}.}%
\label{fig:elliptic_interiority}%
\end{figure}

%% file: 04_numerical_results.tex
%
%
\section{Numerical results}\label{sec:numerical_results}%
In order to assess the proposed algorithm, the present section conducts a two-component series of numerical experiments. %
%
%
Firstly, we investigate various combinations of hypersurfaces and mesh types in subsection~\ref{subsec:results_global}. As a measure for accuracy, we employ the global volume error %
\begin{align}
\error_V\defeq\abs*{1-\frac{\sum_{i=1}^{\Niface}{\abs{\cell*_i\alpha_i}}}{V_\iface}},%
\label{eqn:global_volume_error}%
\end{align}
where $V_\iface=\abs{\set{\vx\in\domain:\lvlset_\iface\fofx\leq0}}$ and $\Niface$ denote the volume enclosed by the hypersurface $\iface$ and the number of intersected cells, respectively. %
The approximation error of the original volume integral in \refeqn{volume_computation_facebased_approx} comprises two distinct sources: (i) the approximation of the hypersurface ($\iface\to\ifaceapprox$ (cf.~subsection~\ref{subsec:metric_projection}) and (ii) the numerical approximation (i.e., quadrature) of the resulting curve integral in \refeqn{quadrature_line_integral}. The latter can be reduced to insignifcance by choosing a sufficiently large order for the employed \textsc{Gauss-Legendre} quadrature. %
\begin{note}[Quadrature]%
In an extensive set of preliminary numerical experiments, we have found that the accuracy of the volume fractions does not profit from increasing $N_{\mathrm{quad}}$ beyond $5$. However, a higher number of quadrature nodes might be needed to accurately approximate general integrals of type $\int_{\iface}{\iprod{\vf}{\nS}\darea}$. %
\end{note}
Hence, the approximation quality of the hypersurface constitutes the limiting factor: %
from the locally quadratic approximation of the hypersurface one can expect third-order convergence with spatial resolution, corresponding to the number of intersected cells $\Niface$, which is not an input parameter. %
The fact that the codimension of $\iface$ with respect to the domain $\domain$ is one implies that $\Niface\propto N_\domain^{\frac{2}{3}}$ or, alternatively, $\sqrt{\Niface}\propto N_\domain^{\frac{1}{3}}$. With $N_\domain^{\frac{1}{3}}$ resembling the equivalent resolution per spatial direction, we choose $\sqrt{\Niface}$ as the corresponding interface resolution. %

%
%
Secondly, note that the meshes under consideration are composed of standard convex polyhedra, which are of high relevance for productive simulations. In order to show the full capability of the proposed algorithm, subsection~\ref{subsec:results_nonconvex} exemplarily investigates a non-convex polyhedron intersected by a family of paraboloids. %

%
%
\subsection{Hypersurfaces}\label{subsec:experiments_hypersurfaces}%
\paragraph{Spheres and ellipsoids}
In this work, we consider a sphere of radius $R_0=\frac{4}{5}$ %
as well as a prolate (semiaxes $\set{\nicefrac{3}{4},\nicefrac{1}{2},\nicefrac{1}{4}}$) and an oblate (semiaxes $\set{\nicefrac{4}{5},\nicefrac{4}{5},\nicefrac{2}{5}}$) ellipsoid, all centered at $\vx_0=\vec{0}$. %
%
%
\paragraph{Perturbed spheres}%
Perturbed spheres can be parametrized in spherical coordinates as %
\begin{align}
\iface=\set{\vx_0+R\fof{\varphi,\theta}\ve_r:\fof{\varphi,\theta}\in\unitsphere}%
\quad\text{with}\quad%
R\fof{\varphi,\theta;\vciface}=\brackets*{\sum\limits_{l=0}^{\Liface}{\sum\limits_{m=-l}^{l}{\ciface{lm}\sphericalharm{}}}}^{\!\frac{1}{3}},%
\end{align}
where the description of the radius $R$ employs tesseral spherical harmonics up to and including order $\Liface\in\setN$. %
The reason for expanding the third power of the radius instead of the radius itself is that the computation of the enclosed volume is considerably simplified, because $V_\iface=\ciface{00}\nicefrac{\sqrt{4\pi}}{3}$ then. %
The $(\Liface+1)^2$ coefficients $\ciface{lm}\sim\mathcal{N}(0,\sigma_0)$ are computed by the method of \citet{AMS_1958_anot}, i.e.
\begin{align}
\ciface{lm}=%
\begin{cases}
\sqrt{4\pi}R_0^3&l=0,\\%
\sqrt{\sigma_0}\sqrt{-2\log\gamma_1}\cos(2\pi\gamma_2)&l>0,%
\end{cases}\qquad\text{with}\qquad\gamma_{1,2}\sim\mathcal{U}(0,1),%
\end{align}
where the uniformly distributed random numbers $\gamma_{1,2}$ are generated by the intrinsic \texttt{fortran} subroutine \texttt{random\_number()}. %
In this work, we consider perturbed spheres with base radius $R_0=\nicefrac{4}{5}$, modes $\Liface\in\set{3,6}$ and variance $\sigma_0=\num{5e-4}$; cf.~\reffig{hypersurface_illustration} for an illustration. %
\begin{figure}[htbp]%
\null\hfill%
\includegraphics[bb=0 0 1190 1100,height=4cm]{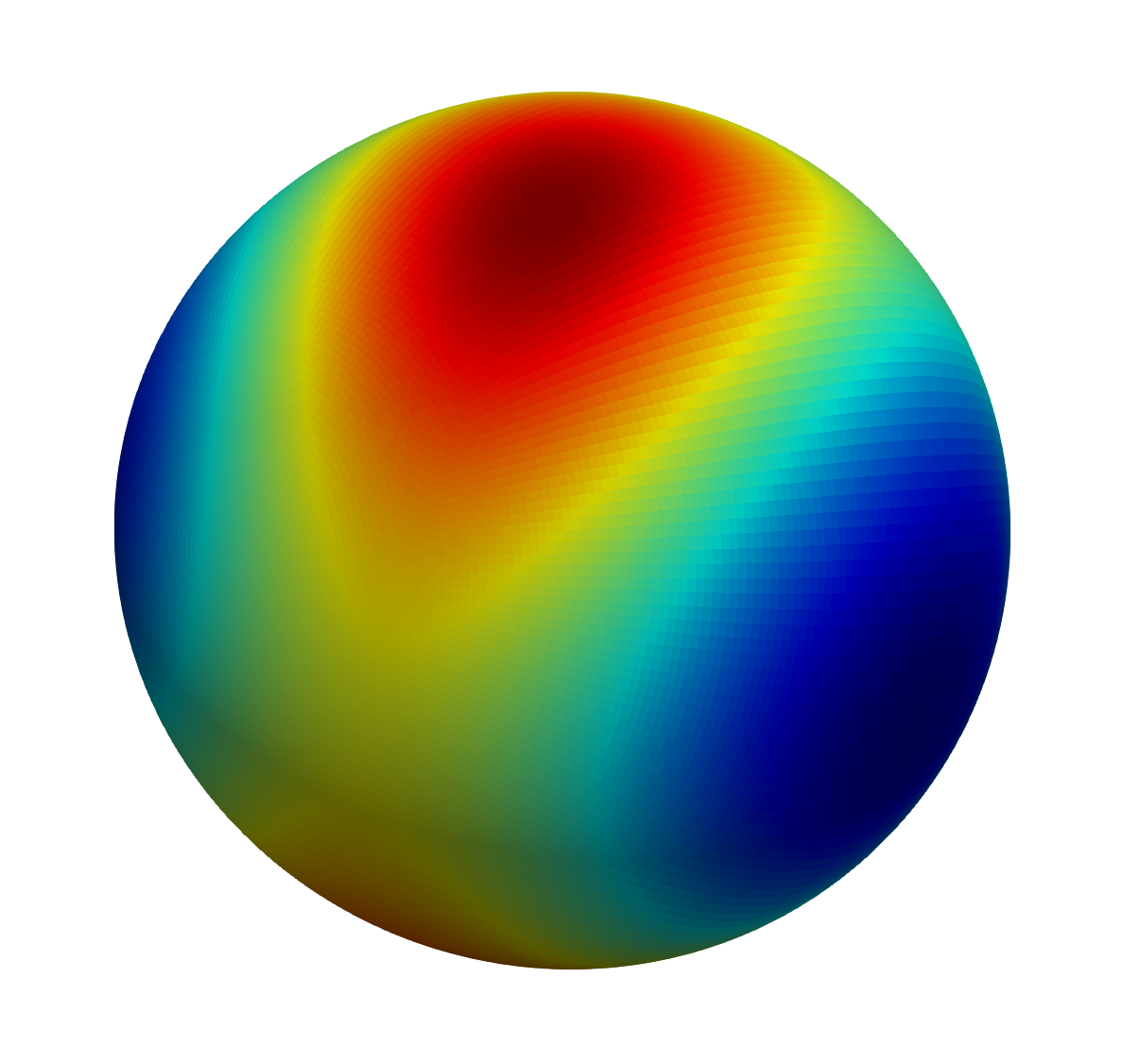}%
\includegraphics[page=1,height=4cm]{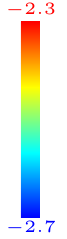}%
\hfill%
\includegraphics[bb=0 0 1190 1100,height=4cm]{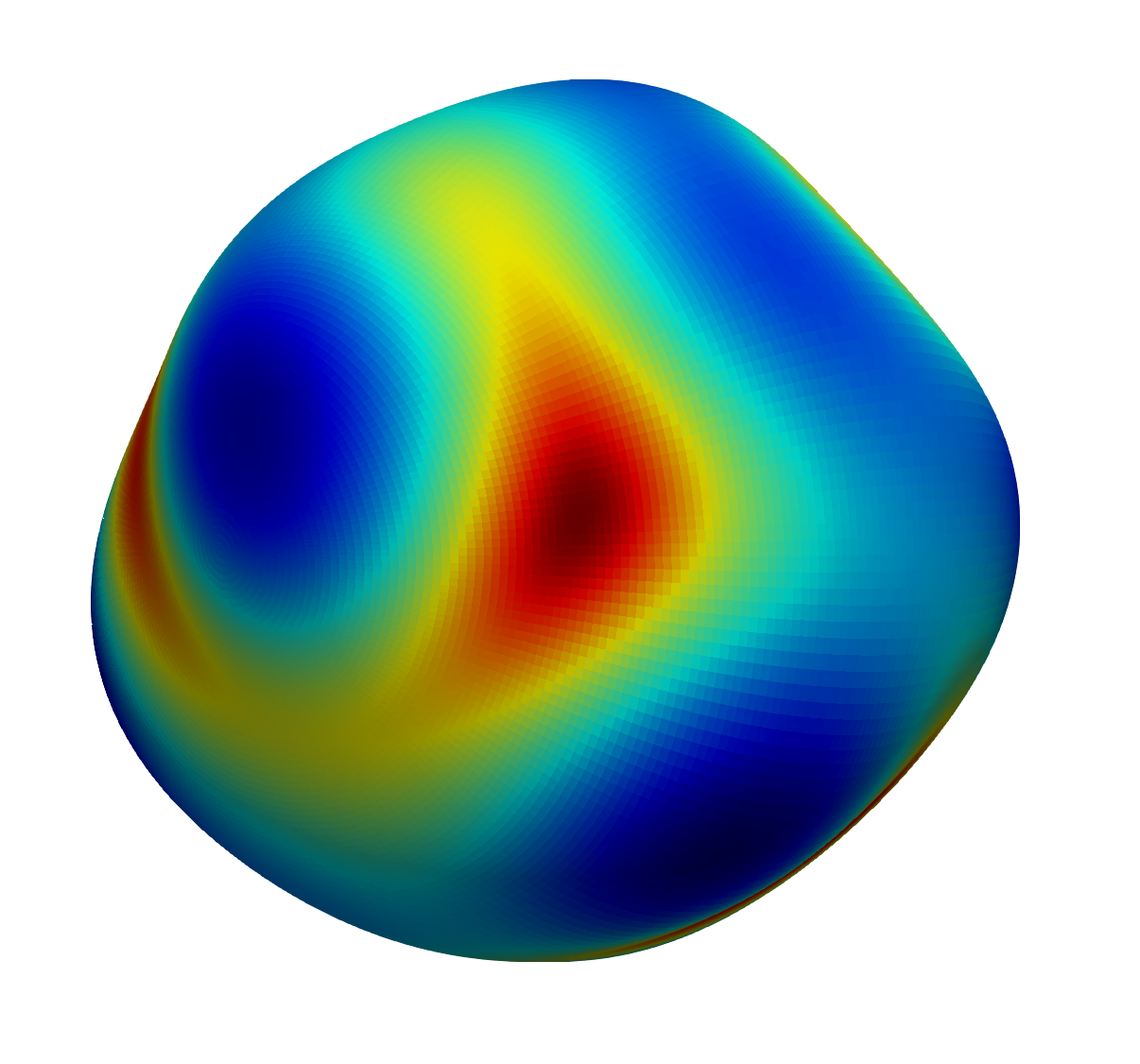}%
\includegraphics[page=2,height=4cm]{colorcurvature}%
\hfill%
\includegraphics[bb=0 0 1190 1100,height=4cm]{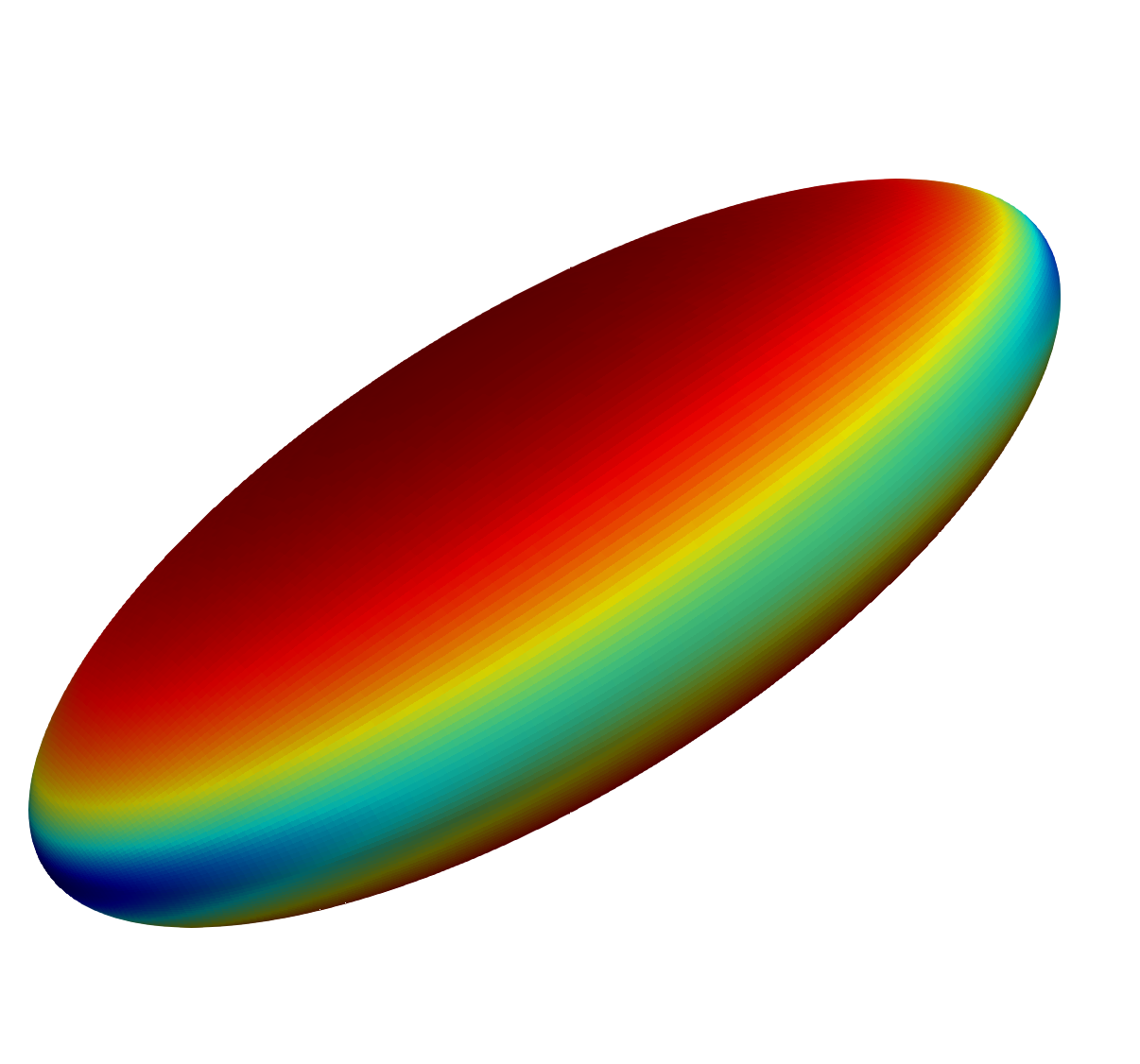}%
\includegraphics[page=3,height=4cm]{colorcurvature}%
\hfill\null%
\caption{Illustration of perturbed spheres ($R_0=\nicefrac{4}{5}$, $\Liface=\protect\set{3,6}$ and $\sigma_0=\num{5e-4}$) and prolate ellipsoid, where the color indicates twice the mean curvature $2\curv_\iface=\curv_1+\curv_2$.}%
\label{fig:hypersurface_illustration}%
\end{figure}
%
%
\subsection{Meshes}\label{subsec:meshes}%
In what follows, we consider the domain $\domain=[-1,1]^3$, which is decomposed into cubes of equal size and tetrahedra. The latter are generated using the library \texttt{gmsh}, introduced in the seminal paper of \citet{IJNME_2009_gmsh}. For the purpose of the present paper, however, we only resort to some of the basic features of \texttt{gmsh}; cf.~\refapp{gmsh} and \reftab{mesh_characteristics} for further details. %

\begin{table}[htbp]%
\caption{Mesh characteristics.}%
\label{tab:mesh_characteristics}%
\null\hfill%
\subtable[Equidistant cube meshes.]{%
\label{tab:cube_mesh_setup}%
\begin{tabular}{c|c|c}%
\textbf{resolution}&\textbf{char. length}&\# of cells\\
$N$&$h=\frac{1}{N}$&$N_\domain=N^3$\\
\hline%
\num{15}&\num{6.66e-02}&\num{3375}\\
\num{20}&\num{5.00e-02}&\num{8000}\\
\num{25}&\num{4.00e-02}&\num{15625}\\
$\vdots$&$\vdots$&$\vdots$\\
\num{70}&\num{1.42e-02}&\num{343000}\\
\end{tabular}%
}%
\hfill%
\subtable[Tetrahedron meshes; cf.~\refapp{gmsh}.]{%
\label{tab:tetrahedron_mesh_setup}%
\begin{tabular}{c|c|c}%
\textbf{resolution}&\textbf{char. length}&\# of cells\\
$N$&$h=\nicefrac{1}{N}$&$N_\domain$\\%
\hline%
\num{10}&$\num{1.00e-1}$&\num{4764}\\
\num{15}&$\num{6.66e-2}$&\num{15266}\\
\num{20}&$\num{5.00e-2}$&\num{33744}\\
\num{25}&$\num{4.00e-2}$&\num{64165}\\
\num{30}&$\num{3.33e-2}$&\num{108582}\\
\num{35}&$\num{2.85e-2}$&\num{171228}
\end{tabular}%
}%
\hfill\null%
\end{table}%
\subsection{Results I -- Meshes with convex cells}\label{subsec:results_global}%
\refFig{volerror_convergence} gathers the volume errors from \refeqn{global_volume_error} obtained for the hypersurfaces and meshes given in subsections~\ref{subsec:experiments_hypersurfaces} and \ref{subsec:meshes}, respectively. %
\begin{figure}[htbp]
\null\hfill%
\subfigure[equidistant cube]{\includegraphics[page=1]{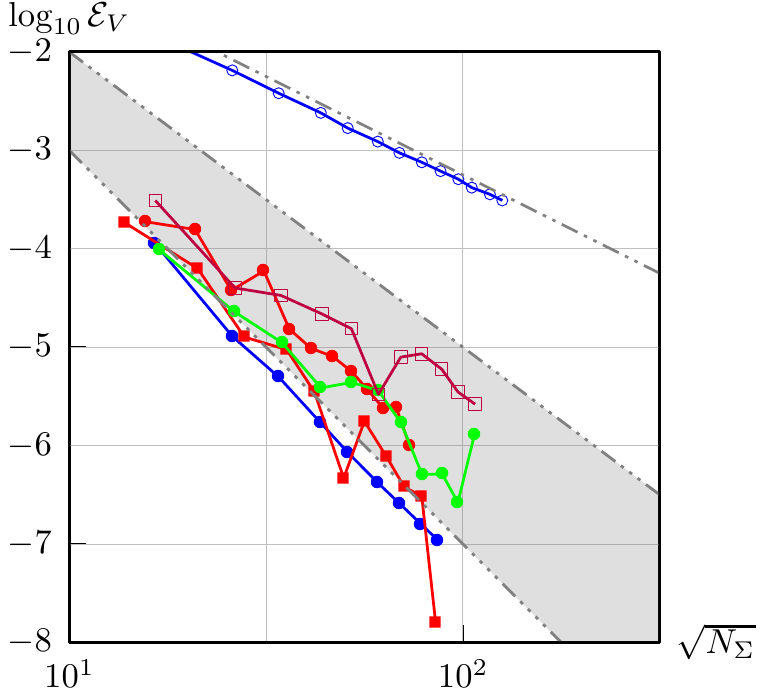}}%
\hfill%
\subfigure[tetrahedron]{\includegraphics[page=2]{volerror_initialization}}%
\hfill\null%
\label{fig:volerror_convergence}%
\caption{Global relative volume error from \protect\refeqn{global_volume_error} as a function of intersected cells $N_\iface$ for different hypersurfaces (centered at $\vx_0=\vec{0}$ in ${\domain=[-1,1]^3}$; cf.~subsection~\ref{subsec:experiments_hypersurfaces} and \protect\reffig{hypersurface_illustration}) and meshes (cf.~\protect\reftab{mesh_characteristics}). %
The number of dots in the dashed lines corresponds to the order of convergence. For comparison, we have added the results for a sphere with linear hypersurface approximation (\textcolor{blue}{$\circ$}).}%
\end{figure}%

The main observations can be summarized as follows: %
\begin{enumerate}
\item As expected, the global relative volume error exhibits at least third-order convergence with spatial resolution for all combinations of hypersurfaces and meshes under consideration. %
\item For spheres, one obtains fourth-order convergence for both tetrahedral and hexahedral meshes. The rationale behind this phenomenon emerges from considering the approximation quality of the height function in \refeqn{hypersurface_approximation_quadratic}: for general hypersurfaces, the quadratic approximation in a tangential coordinate system exhibits third-order. Due to the symmetry, however, the general remainder $\hot{\norm{\vt}^3}$ effectively becomes $\hot{\norm{\vt}^4}$, which directly translates to an increased order of convergence. %
\item While the error obviously decreases with increasing spatial resolution, it is virtually independent of the underlying mesh, indicating the robustness of the proposed method. %
\item For comparison, \reffig{volerror_convergence} also contains the volume errors obtained from linear hypersurface approximation ($\curv_i\defeq0$). For both cube and tetrahedral meshes, the error differs between two and four orders of magnitude, which is in accordance with the findings of \citet{JCP_2019_haco}. Note that, while we only show the results for the sphere, they can be considered prototypical for the other hypersurfaces under consideration. However, owing to the approximation quality discussed above, one obtains a reduced difference (two to three) in the order of magnitude. %
\item For cube meshes, the results virtually coincide with those of \citet{JCP_2019_haco}. Due to the strong similarity in concept, this is to be expected. However, recall that the proposed method is applicable to unstructured meshes composed of arbitrary polyhedra, whereas the original algorithm in \cite{JCP_2019_haco} is restricted to (i) convex polyhedral cells enclosing (ii) simply connected hypersurface patches. %
\end{enumerate}%
\paragraph{Computational time}%
In addition to the accuracy of the volume approximation, \reftab{component_timing} assesses the performance of the proposed algorithm in terms the componentwise share of computational time. %
\begin{table}[htbp]%
\centering%
\renewcommand{\arraystretch}{1.25}%
\caption{Computational time of sub-algorithms (averaged over intersected cells) in \% of the total computational time (I) and in multiples of the reference time $T_\mathrm{EVD}$ (II; see~\refapp{reference_time}).}%
\label{tab:component_timing}%
\begin{tabular}{|c|c|c|c|}%
\hline
\textbf{sub-algorithm}&\textbf{info}&I&II\\%
\hline
{hypersurface approximation}&\ref{subsec:approximation_hypersurfaces},\ref{subsec:metric_projection}&25.31\%&\num{1.17E-03}\\%
{face intersection}&\ref{subsec:computation_edge_intersections}&12.51\%&\num{5.79E-04}\\%
{principal transformation}&\ref{subsec:transformation_principal_coordinates}&32.93\%&\num{1.50E-03}\\%
{interiority check (elliptic)}&\reffig{elliptic_interiority}&3.54\%&\num{1.66E-04}\\%
{reconstrution of} $\partial\ifaceapprox_k$&\ref{subsec:topological_connectivity}, \ref{subsec:computation_immersed_area}&19.41\%&\num{9.04E-04}\\%
{quadrature (evaluation)}&\refeqn{quadrature_line_integral}&6.30\%&\num{2.95E-04}\\%
\hline
\end{tabular}%
\end{table}%

\subsection{Results II -- Single non-convex polyhedron}\label{subsec:results_nonconvex}%
The previous subsection~\ref{subsec:results_global} was devoted to the investigation of general hypersurfaces intersecting hexahedral and tetrahedral (i.e., convex) meshes, highlighting the influence of hypersurface approximation; cf.~subsection~\ref{subsec:metric_projection}. In addition, the present section focusses on the volume computation for a given family of paraboloids intersecting a single non-convex polyhedron; cf.~\refapp{nonconvex_polyhedron} for details. %
\paragraph{A family of paraboloids}%
For a given base point $\xref$, base system $\set{\vn_0,\vtau}$ and curvature tensor $\curvtensor$, extending \refeqn{ifaceapprox_map} by a shift in the direction of the base normal $\vn_0$ yields a family of paraboloids, namely %
\begin{align}
\set{\ifaceapprox}_{s\in[s_{\mathrm{min}},s_{\mathrm{max}}]}=\set{\vf_\ifaceapprox\fof{\vt;s}:(\vt,s)\in\setR^2\times[s_{\mathrm{min}},s_{\mathrm{max}}]}%
\quad\text{with}\quad%
\vf_\ifaceapprox\fof{\vt;s}=\xref+\vtau\vt+\brackets{\happrox\fof{\vt}+s}\vn_0,%
\tag{\ref*{eqn:ifaceapprox_map}$^\prime$}%
\label{eqn:shifted_paraboloid}%
\end{align}
where the associated level-set analogously extends \refeqn{ifaceapprox_lvlset}. %
As parameters of the paraboloid, we choose %
\begin{align}%
\xref=\frac{\brackets*[s]{1,1,1}\transpose}{2},\quad%
\vn_0=\frac{\brackets*[s]{4,-7,2}\transpose}{\sqrt{69}},\quad%
\vtau_1=\frac{\brackets*[s]{-8,14,65}\transpose}{\sqrt{4485}},\quad
\vtau_2=\frac{\brackets*[s]{-7,-4,0}\transpose}{\sqrt{65}},\quad%
\curv_1=-\frac{19}{4}\quad\text{and}\quad%
\curv_2=0.%
\label{eqn:paraboiloid_parameters}%
\end{align}
The left panel in \reffig{illustration_graphbase} illustrates the intersection with a unit cube for $s=0$. %
\begin{remark}[Choice of parameters]%
Recall from \reftab{possible_curve_class} that there are four classes of boundary curves: hyperbolic, elliptic, parabolic and linear. %
The first two require (i) non-zero \textsc{Gaussian} curvature as well as (ii) non-orthogonality of the containing face $\vfface_k$ and the base plane of the paraboloid, i.e.~$\iprod{\vn^\vfface_k}{\vn_0}\neq0$. Complementary, parabolic curve segments may emerge if either (i) one of the principal curvatures $\curv_i$ is zero for $\iprod{\vn^\vfface_k}{\vn_0}\neq0$ or (ii) $\iprod{\vn^\vfface_k}{\vn_0}=0$ for arbitrary values of $\curv_i$.
Hence, due to the choice of the hypersurfaces in subsection~\ref{subsec:experiments_hypersurfaces}, in statistical terms, one cannot expect to encounter parabolic or linear boundary curve segments. %
Therefore, the present subsection aims at examining parabolic and linear curve segments by purposely setting one of the principal curvatures, say $\curv_2$, to zero. %
\end{remark}%
The boundaries of the shift interval are chosen to ensure that $\polyvof\fof{s_{\mathrm{min}}}=0$ and $\polyvof\fof{s_{\mathrm{max}}}=1$, i.e.~such that the volume fraction
\begin{align}
\polyvof\fof{s}=\abs{\cell*}^{-1}\abs{\set{\vx\in\cell*:\lvlset_\ifaceapprox\fof{\vx;s}\leq0}}%
\label{eqn:nonconvex_volume_fraction}%
\end{align}
traverses all possible values. Here, let $s_{\mathrm{min}}\defeq-1$ and $s_{\mathrm{max}}\defeq\nicefrac{3}{2}$. %
It is worth noting that, for non-degenerate\footnote{\citet[Section~2]{JCP_2021_fbip} consider the regularity of planar $\ifaceapprox$ for intersection with convex and non-convex polyhedra.} paraboloids ($\curv_i\neq0$), the function $\polyvof:\setR\mapsto[0,1]$ is strictly monotonous and continuously differentiable. The regularity in the degenerate case (at least one trivial principal curvature) depends on the topological properties of the polyhedron as well as the paraboloid parameters. %
\begin{remark}[Partial derivative]%
After replacing $\iface$ with $\ifaceapprox$ in \refeqn{volume_computation} and applying the \textsc{Reynolds} transport theorem, it is easy to show that %
$\abs{\cell*}\partial_s\polyvof\fof{s}=\int_{\ifaceapprox}{\iprod{\partial_s\vf_\ifaceapprox}{\vn_\ifaceapprox}\darea}=\abs{\support_\ifaceapprox}\geq0$, %
i.e.~the derivative of the volume fraction with respect to the base normal shift parameter corresponds to the area of the graph base $\support_\ifaceapprox$. For planar paraboloids, one trivially obtains $\abs{\support_\ifaceapprox}=\abs{\ifaceapprox\cap\cell*}$, which can be exploited, e.g., for efficient PLIC interface positioning schemes \cite{XXX_2020_ivof,JCP_2021_fbip,JCP_2021_etmp}. 
\end{remark}
\refFig{nonconvex_volume_fraction} depicts the volume fraction $\polyvof$ and its derivative with respect to the shift parameter $s$ as a function of the latter, where \reffig{nonconvex_volume_fraction_illustration} illustrates some of the intersections. %
\begin{figure}[htbp]%
\null\hfill%
\includegraphics{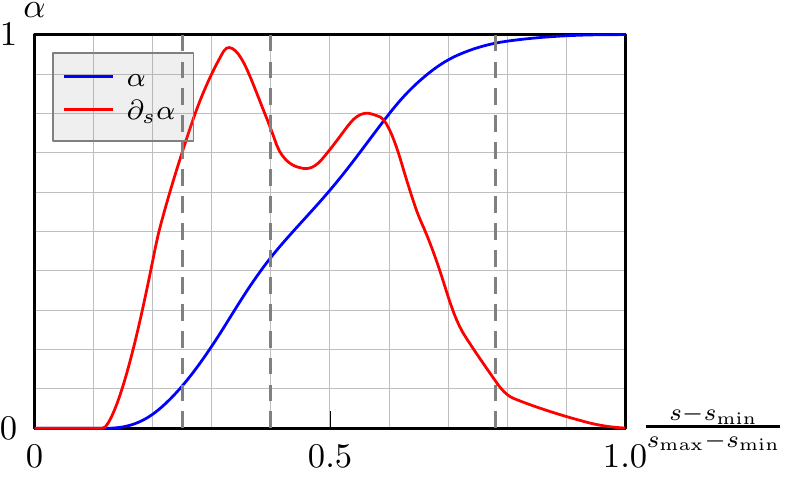}%
\hfill\null%
\caption{Volume fractions $\polyvof$ from \protect\refeqn{nonconvex_volume_fraction} induced by the family of paraboloids given in \protect\refeqn{shifted_paraboloid} as a function of the shift $s$ with ${[s_{\mathrm{min}},s_{\mathrm{max}}]=[-1,\nicefrac{3}{2}]}$.}%
\label{fig:nonconvex_volume_fraction}%
\end{figure}%

\begin{figure}[htbp]%
\null\hfill%
\includegraphics[page=1]{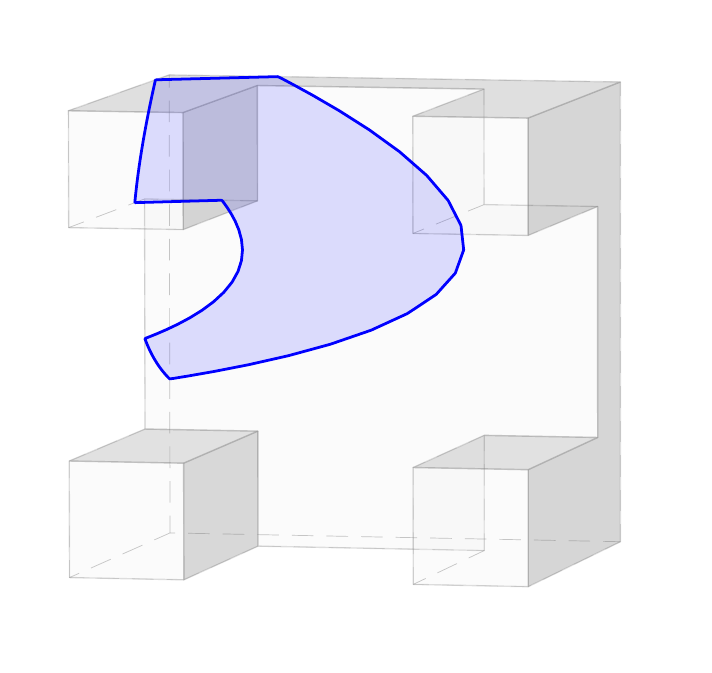}%
\hfill%
\includegraphics[page=1]{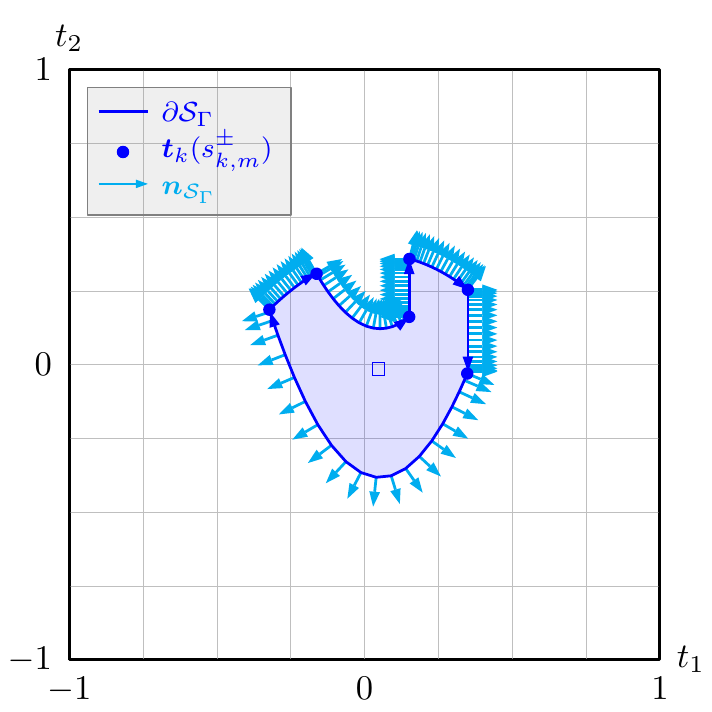}%
\hfill\null%
\\%
\null\hfill%
\includegraphics[page=2]{nonconvex_example}%
\hfill%
\includegraphics[page=2]{nonconvex_boundarycurve_paradom}%
\hfill\null%
\\%
\null\hfill%
\includegraphics[page=3]{nonconvex_example}%
\hfill%
\includegraphics[page=3]{nonconvex_boundarycurve_paradom}%
\hfill\null%
\caption{Polyhedron "table" intersected by shifted paraboloid (blue shade: $\ifaceapprox\cap\protect\cell*$) given in \protect\refeqn{shifted_paraboloid} with associated graph base $\support_\ifaceapprox$ at normalized positions $\nicefrac{1}{4}$, $\nicefrac{2}{5}$ and $\nicefrac{39}{50}$ (corresponding to the vertical dashed lines in \protect\reffig{nonconvex_volume_fraction}) in interval ${[s_{\mathrm{min}},s_{\mathrm{max}}]=[-1,\nicefrac{3}{2}]}$.}%
\label{fig:nonconvex_volume_fraction_illustration}%
\end{figure}

%% file: 05_conclusion.tex
\section{Conclusion \& Outlook}\label{sec:conclusion}%
We have introduced an algorithm for the computation of volumes induced by an intersection of a paraboloid with an arbitrary polyhedron, where the paraboloid parameters are obtained from a locally quadratic approximation of a given hypersurface. %
The recursive application of the \textsc{Gaussian} divergence theorem in its respectively appropriate form allows for a highly efficient face-based computation of the volume of the truncated polyhedron, implying that no connectivity information has to be established at runtime. Furthermore, the face-based character renders the presented approach most suitable for parallel computations on unstructured meshes. %
A classification of the boundary curve segments $\partial\ifaceapprox_k$ associated to the polyhedron faces $\vfface_k$ allows for their explicit parametrization, which has been shown to be favourable for the computation of quadrature nodes and weights. This in turn strongly facilitates the approximation of the associated curve integrals. %
We have conducted a twofold assessment of the proposed algorithm: %
firstly, by examining convex meshes with different hypersurfaces, resemebling a commonly encountered task for obtaining the initial configuration of two-phase flow simulations. For each cell, the paraboloid parameters are obtained from a locally quadratic approximation based on the level-set of the original hypersurface. %
The global volume errors show the expected third- to fourth-order convergence with spatial resolution, along with an error reduction of about 2 orders of magnitude in comparison to linear approximation. %
Secondly, by intersecting a paramterized family of paraboloids with an exemplary non-convex polyhedron, which serves to illustrate the capability of the proposed algorithm. %

Altogether, we draw the following conclusions: %
\begin{enumerate}%
\item The recursive application of the \textsc{Gaussian} divergence theorem in appropriate form allows for an efficient computation of the volume of a truncated \textbf{arbitrary polyhedron}. This face-based decomposition allows to \textbf{avoid extracting topological connectivity}, which is advantageous in terms of implementation complexity, computational effort and parallelization. %
\item The quadrature nodes and weights can easily be employed to evaluate general integrals of type $\int_{\ifaceapprox}{\iprod{\vf}{\vn_\ifaceapprox}\darea}$ for integrands $\vf$ which are polynomial in the spatial coordinate $\vx$. Note that, e.g., the partial derivatives of the volume with respect to the paraboloid parameters can be written in this form. In fact, the present algorithm constitutes an important building block for a generalization of the parabolic reconstruction of interfaces from volume fractions, originally proposed for structured hexahedral grids by \citet{JCP_2002_prost}. %
\end{enumerate}%

%% file: 99_appendix.tex
\begin{appendix}
%
%
\section{A machine-independent reference for computational time}\label{app:reference_time}%
We seek to establish a referential measure for computational time which is both easily reproducable and obtainable in most technically relevant programming languages. %
Computing the eigenvalues of a real non-symmetric matrix constitutes a frequent task in many fields of physics, where open-source libraries such as \texttt{LAPACK} contain highly efficient implementations; cf.~\citet{lapack99}. %
For the purpose of the present study, we compute the eigenvalues and right eigenvectors of the matrix %
\begin{align}
\tensor[2]{M}=%
\brackets*[s]{%
\begin{matrix}%
5&3&0&-100&2\\%
0 &6 &0&11 & 0\\%
1 &1 &0&-2 & 0\\%
99&-4&0&-10& 2\\%
1 &0 &7&  6&-4%
\end{matrix}}%
\end{align}
using the \href{http://www.netlib.org/lapack/explore-html/d9/d8e/group__double_g_eeigen_ga66e19253344358f5dee1e60502b9e96f.html}{routine \texttt{DGGEV}}\footnote{Note that the input-parameter \texttt{lwork} was determined by the recommended query run.} via %
\begin{center}
\verb+DGEEV('N','V',5,+$\tensor[2]{M}$\verb+,5,wr,wi,vr,5,vr,5,work,lwork,info)+.%
\end{center}
For reasons of robustness, we consider the total execution time, say $T_\mathrm{EVD}$, of $\num{e3}$ calls. %
%
%
\section{Mesh generation with \texttt{gmsh}}\label{app:gmsh}%
The tetrahedral meshes used in section~\ref{sec:numerical_results} were generated with \href{https://gmsh.info}{\texttt{gmsh} 4.7.1}. For $h=\frac{1}{N}$, cf.~\reftab{tetrahedron_mesh_setup}, the file \texttt{mesh.geo} gathers the relevant information. %
\lstset{language=C++,commentstyle=\color{blue!50},frame=single,caption={Example geometry file ($N=20$)}}%
\begin{lstlisting}
    // add points
    Point(1)={0,0,0,h};Point(2)={1,0,0,h};Point(3)={1,1,0,h};Point(4)={0,1,0,h};
    Point(5)={0,0,1,h};Point(6)={1,0,1,h};Point(7)={1,1,1,h};Point(8)={0,1,1,h};
    // add lines
    Line(1) = {1, 2};Line(2) = {2, 3};Line(3) = {3, 4};Line(4) = {4, 1};
    Line(5) = {1, 5};Line(6) = {2, 6};Line(7) = {3, 7};Line(8) = {4, 8};
    Line(9) = {5, 6};Line(10) = {6, 7};Line(11) = {7, 8};Line(12) = {8, 5};
    // loops and surfaces
    Line Loop(13) = {5,-12,-8,4};Line Loop(14) = {2,7,-10,-6};
    Line Loop(15) = {-4,-3,-2,-1};Line Loop(16) = {9,10,11,12};
    Line Loop(17) = {1,6,-9,-5};Line Loop(18) = {3,8,-11,-7};
    Plane Surface(1) = {13};Plane Surface(2) = {14};Plane Surface(3) = {15};
    Plane Surface(4) = {16};Plane Surface(5) = {17};Plane Surface(6) = {18};
    Surface Loop(1) = {6, 3, 1, 5, 2, 4};Volume(1) = {1};
\end{lstlisting}
With the above geometry file, the mesh is generated by invoking %
\begin{center}
\verb+gmsh -refine -smooth 100 -optimize_netgen -save -3 -format vtk -o mesh.vtk mesh.geo+
\end{center}
\section{A non-convex polyhedron}\label{app:nonconvex_polyhedron}
\begin{figure}[htbp]
\null\hfill%
\includegraphics{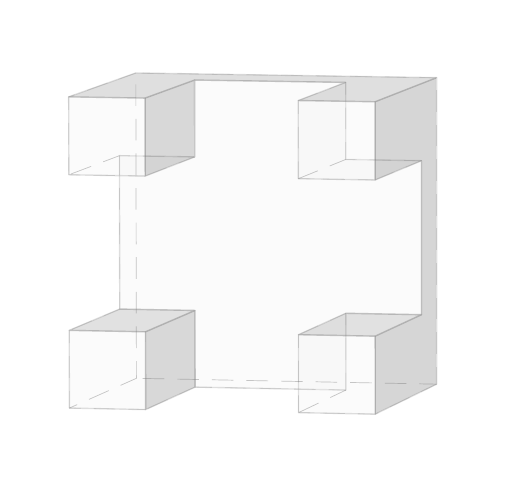}%
\hfill\null%
\caption{Non-convex polyhedron \textit{table}.}
\label{fig:table_illustration}%
\end{figure}

We consider a table-shaped polyhedron composed of a cuboid "plate" of size $1\times1\times a$ and four cuboid "legs" of size $a\times a\times1-a$ with $a=\nicefrac{1}{4}$, as illustrated in \reffig{table_illustration}. Note that the polyhedron contains non-convex faces. %
\end{appendix}